\documentclass[aop,MSNbibl,seceqn,citesort,dvips]{arximspdf}
\usepackage{mathrsfs}

\doi{10.1214/11-AOP676}
\volume{40}
\issue{5}
\pubyear{2012}
\firstpage{1945}
\lastpage{1979}

\makeatletter

\newtheorem{theorem}{Theorem}[section]
\newtheorem{prop}[theorem]{Proposition}
\newtheorem{cor}[theorem]{Corollary}
\newtheorem{lem}[theorem]{Lemma}
\newproclaim{remark}[theorem]{Remark}
\newproclaim{assumption}{Assumption}

\newcommand{\DD}{/\!\!/}

\newcommand{\PP}{\mathbb{P}}
\newcommand{\E}{\mathbb{E}}
\newcommand{\R}{\mathbb{R}}
\newcommand{\N}{\mathbb{N}}
\newcommand{\sC}{\mathscr{C}}
\newcommand{\sD}{\mathscr{D}}
\newcommand{\sF}{\mathscr{F}}
\newcommand{\sL}{\mathscr{L}}

\newcommand{\gm}{\gamma}
\newcommand{\ep}{\varepsilon}
\newcommand{\zt}{\zeta}
\newcommand{\h}{\eta}
\newcommand{\kp}{\kappa}
\newcommand{\ro}{\rho}
\newcommand{\sg}{\sigma}

\newcommand{\Dl}{\Delta}
\newcommand{\Lm}{\Lambda}
\newcommand{\Ph}{\Phi}

\newcommand{\wg}{\wedge}
\newcommand{\nab}{\nabla}
\newcommand{\Ric}{\operatorname{Ric}}
\newcommand{\Cut}{\operatorname{Cut}}
\newcommand{\Cutst}{\Cut_{\mathrm{ST}}}

\makeatother

\begin{document}
\begin{frontmatter}

\title{Convergence of time-inhomogeneous geodesic~random walks
and its application to~coupling methods}
\runtitle{Time-inhomogeneous geodesic random walks}

\begin{aug}
\author[A]{\fnms{Kazumasa} \snm{Kuwada}\corref{}\thanksref{t1}\ead[label=e1]{kuwada.kazumasa@ocha.ac.jp}}
\runauthor{K. Kuwada}
\affiliation{Ochanomizu University}
\address[A]{Graduate School of Humanities and Sciences\\
Ochanomizu University\\
Ohtsuka 2-1-1, Bunkyo-ku\\
Tokyo 112-8610\\
Japan\\
\printead{e1}} 
\end{aug}

\thankstext{t1}{Supported in part
by the JSPS fellowship for research abroad and by KAKENHI (22740083).}

\received{\smonth{7} \syear{2010}}
\revised{\smonth{3} \syear{2011}}

%
\begin{abstract}
We study an approximation by time-discretized geodesic random walks of
a diffusion process associated with a family of time-dependent metrics
on manifolds. The condition we assume on the metrics is a natural
time-inhomogeneous extension of lower Ricci curvature bounds. In
particular, it includes the case of backward Ricci flow, and no further
a priori curvature bound is required. As an application, we construct a
coupling by reflection which yields a nice estimate of coupling time,
and hence a gradient estimate for the associated semigroups.
\end{abstract}

%

\end{frontmatter}

\section{Introduction}
\label{sec1}


It has been well known that there is a strong connection
between behavior of heat distributions or Brownian motions
and geometry of their underlying space.
Even on time-inhomogeneous spaces such as Ricci flow,
this guiding principle has been confirmed
through recent developments
(see
\cite
{Arn-Coul-Thalhoriz,CoulgtBM,K-Phili2,K-Phili,McC-ToppWass-RF,ToppLopt,ZhangQSbook}
and references therein).
Some of them~\cite{Arn-Coul-Thalhoriz,K-Phili2}
are based on coupling methods of stochastic processes.
Given two stochastic processes $Y_1 (t)$ and $Y_2 (t)$
on a state space $M$,
a coupling $\mathbf{X} (t) = ( X_1 (t) , X_2 (t) )$
of $Y_1 (t)$ and $Y_2 (t)$
is a stochastic process on $M \times M$ such that
$X_i$ has the same law as $Y_i$ for $i=1,2$.
By constructing a suitable coupling which reflects
the geometry of the underlying structure,
one can obtain various estimates
for heat kernels, harmonic maps, eigenvalues
etc. under natural geometric assumptions
(see~\cite{Hsu,Kendsurvey,Wangbook05}, e.g.).
Since coupling of random variables
provides a coupling of their distributions,
coupling methods are naturally connected
with the theory of optimal transportation, which
are used in some of aforementioned results
\cite{McC-ToppWass-RF,ToppLopt}.
With further studies in this direction in mind,
here we consider an approximation of diffusion processes
associated with a~family of time-dependent metrics
by so-called geodesic random walks.
Generally speaking, one of the major reasons
that we establish approximation
is to overcome technical difficulties
in studying the object in the limit.
This is also our case, and we will use the approximation
in order to study a~coupling of diffusion processes.

Let $M$ be a smooth manifold with a family of
complete Riemannian metrics $g(t)$
indexed by $t \in[T_1 , T_2]$.
By $( X (t) )_{t \in[ T_1 ,T_2 ] }$,
we denote the $g(t)$-Brownian motion.
It means that
$X (t)$ is a time-inhomogeneous
diffusion process on $M$
associated with $\Delta_{g(t)} /2$,
where $\Delta_{g(t)}$ is the Laplacian
with respect to $g(t)$
(see~\cite{CoulgtBM} for a construction of $g(t)$-Brownian motion).
A geodesic random walk $\tilde{X}^\alpha$ on $M$
with a parameter $\alpha$
is a discrete time Markov chain
whose one-step variation is given as follows:
Given a position $x$ at some time~$t$,
consider a random vector in $T_x M$.
We map it to $M$ by $g(t)$-exponential map
to determine the next position.
Here the parameter $\alpha$ is implemented
as a (diffusive) scaling
on time step and on the length of the random vector
in $T_x M$; see Section~\ref{sec3} for more details.
In this paper, we consider only the case that
all the random vectors in tangent space is specified
to the one having a uniform distribution
on a $g(t)$-ball whose radius is comparable to~$\alpha$.
A~simplified version of our main theorem,
the convergence of geodesic random walks,
is stated as follows;
see Theorem~\ref{thIP} and Section~\ref{sec3}
for a more precise and general statement:
\begin{theorem} \label{thmain00}
Suppose
%
\begin{equation} \label{eqc-b}
\partial_t g(t) \le\Ric_{g(t)}
\end{equation}
holds.
Then a continuous time interpolation of $\tilde{X}^\alpha$
converges in law to $X$ as $\alpha\to0$.
\end{theorem}

As we will see in the sequel,
there are several technical difficulties
arising from the time-dependency on the metric.
Nevertheless, the assumption of
the full statement in Theorem~\ref{thIP}
is much weaker in some respect than that
in the classical time-homogeneous case.
Thus this assertion itself would be of interest,
independently of its application to coupling methods.

In the time-homogeneous case,
the convergence in law of
scaled geodesic random walks
to the Brownian motion
is used to study a coupling of Brownian motions
$( X_1 (t) , X_2 (t) )$ by reflection;
see~\cite{K8,Renespoly}.
A coupling of this kind provides us a useful control
of the coupling time $\tau^*$,
the first time when $X_1$ and $X_2$ meet.
Even in our time-inhomogeneous case,
Theorem~\ref{thmain00} carries
the same estimate in almost the same way.
A simplified version of this assertion
is as follows;
for the complete statement of our main theorem,
see Theorem~\ref{thmain0}.
%
\begin{theorem} \label{thmain}
Suppose~(\ref{eqc-b}).
Then, for each $x_1 , x_2 \in M$,
there exists a coupling
$\mathbf{X} (t) := ( X_1 (t) , X_2 (t) )$
of two $g(t)$-Brownian motions
starting at $(x_1 , x_2 )$
satisfying
%
\begin{equation} \label{eqdom}
\PP[ \tau^* > t ]
\le
\PP\biggl[ \inf_{T_1 \le s \le t} B (s) > - \frac{d_{g(T_1)} ( x_1 ,
x_2 )}{2} \biggr]
\end{equation}
for each~$t$, where $d_{g(T_1)}$ is the distance
function on $M$ with respect to $g(T_1)$, and
$B (t)$ is a one-dimensional standard Brownian motion
starting at the time $T_1$.
\end{theorem}

Similarly to the time-homogeneous case,
Theorem~\ref{thmain} yields
a gradient estimate of the heat semigroup,
which tells us a quantitative estimate
on the smoothing effect of the heat semigroup;
see Corollary~\ref{corsF}.
In addition, we can apply our method
to construct different kinds of couplings.
As one of them, coupling by space--time parallel transport
is studied in~\cite{K-Phili2} by using Theorem~\ref{thmain00},
and it sharpens the result by Topping~\cite{ToppLopt},
concerning the monotonicity of a transportation cost
between the heat distributions whose cost is measured
by Perelman's $\sL$-distance.

Condition~(\ref{eqc-b}) is essentially
the same as backward super Ricci flow in~\cite{McC-ToppWass-RF}; our
condition is slightly different in constant
since our $g(t)$-Brownian motion,
and hence our heat equation
corresponds to $\Dl_{g(t)} / 2$ instead of $\Dl_{g(t)}$.
Obviously,~(\ref{eqc-b}) is satisfied if $g(t)$ evolves
according to the backward Ricci flow
$\partial_t g(t) = \Ric_{g(t)}$.
From a different point of view, condition~(\ref{eqc-b}) can be
interpreted as
a time-inhomogeneous analog of nonnegative Ricci curvature
since $\partial_t g(t)$ vanishes if $g(t)$ is independent of~$t$.
Along this viewpoint, we can consider
a time-inhomogeneous analog of
more general lower Ricci curvature bounds,
and we can obtain the conclusion under
such generalized conditions in the sequel; see Assumption \ref
{assnon-explosion} and
(\ref{eqc-b0}); cf. Remark~\ref{remO-U}.
It should be remarked that,
even in those cases, no uniform lower bound of $\Ric_{g(t)}$
only in terms of $g(t)$ without time derivative
is assumed.
In particular,
no bounds of $g(t)$ or $g(t)$-curvature tensor
being uniform in time are required.
Since a~Ricci flow will produce a singularity in a~finite time,
a~time-uniform bound on $g(t)$ or $\Ric_{g(t)}$
seems to be restrictive.
It might be possible to simplify the proof by
supposing additional assumptions
involving a time-uniform estimate; however, this is out of the scope of
this paper.

In our argument, the distance function $d_{g(t)}$
with respect to the time-dependent metric $g(t)$
plays a prominent role.
The first variation in~$t$ of $d_{g(t)}$ is described
in terms of $\partial_t g(t)$, and the second variation of
$d_{g(t)}$ in space variables involves a notion of curvature.
Both of these variations appear in the bounded variation part of
the radial process $d_{g(t)} ( o , X(t) )$ of
the $g(t)$-Brownian motion via the It\^o formula.
Thus a relation between
$\partial_t g(t)$ and $\Ric_{g(t)}$ [e.g.,~(\ref{eqc-b})], produces a
nice control
of the radial process.
Although we will work on geodesic random walks instead of
the $g(t)$-Brownian motion itself,
such an observation is still efficient.\vadjust{\goodbreak}

In the time-homogeneous case,
the convergence of scaled geodesic random walks
is first studied by J{\o}rgensen~\cite{Jorg}
by using the convergence theory of semigroups; see \cite
{Blum,PingRW}, also.
However, in our framework,
it is not clear whether we can apply a similar technique
since the base measure, the Riemannian volume,
depends on time, and hence we cannot expect that
it makes the heat semigroup invariant.
To avoid such a technical difficulty,
we use the uniqueness of the martingale problem instead
for identifying the limit.
Another difficulty arises from
the lack of time-uniform bounds of Riemannian metrics.
It prevents us to expect a \textit{global} comparison of
geometric structures,
such as $d_{g(t)}$, between different times.
Thus we will make some efforts for localizing the problem
by giving a uniform estimate of the first exit time of~$X^\alpha$
from a large ball centered at a reference point.
Note that our assumption admits
lower unbounded Ricci curvatures
even in the case $\partial_t g (t) \equiv0$
(see Assumption~\ref{assnon-explosion}).
Thus our assumption on the geometry of
the underlying space is weaker than that in~\cite{Jorg}
(by considering Riemannian manifolds
with lower unbounded curvature,
we can easily find an example
which does not satisfy the assumption in
\cite{Jorg}).
On the other hand, the assumption
on the driving noises of the geodesic random walk
in~\cite{Jorg}
is more general than our specified one.
Though it might be possible
to take a more general noise under our assumption,
our result already works well
for applying the approximation to coupling methods.
As a related work,
the theory of time-dependent Dirichlet forms
has been developed for studying
the time-inhomogeneous Markov processes
in the literature~\cite{OshimatdepDF}; see also~\cite{StangDF}.
Unfortunately, because of above-mentioned difficulties,
our framework does not fall into the scope of
those theories at this moment.
It might be an important problem to extend those theories
so that they includes our case.

The organization of this paper is as follows.
In the rest of this section after this paragraph,
we review existing approaches
on the construction of couplings.
By comparing those approaches with ours,
we try to explain the reason
why we choose our approach
for constructing a coupling by reflection.
In the next section,
we show basic properties
of a family of Riemannian manifolds
$( ( M , g(t) ) )_{t}$.
In particular, we prove that
Riemannian metrics $( g(t) )_t$ are
\textit{locally} comparable with each other.
It will be used to give a uniform control of
several error terms which appear
as a result of our discrete approximation.
In Section~\ref{sec3}, we will study geodesic random walks
in our time-inhomogeneous framework.
There we introduce them and
prove the convergence in law
to a diffusion process.
After a small discussion
at the beginning of the section,
the proof is divided into two main parts.
In the first part,
we will give a uniform estimate
for the exit time of geodesic random walks
from a big compact set.
Our assumption here is almost the same as
in~\cite{K-Phili} where nonexplosion
of the diffusion process is studied; see Remark~\ref{remtight} for
more details.
In the second part, we prove the tightness of geodesic random walks
on the basis of the result in the first part.
In Section~\ref{sec4}, we will construct a coupling by reflection
and show an estimate of coupling time,
which completes the proof of Theorem~\ref{thmain} as a special case.
In Section~\ref{sec5}, we will give a short remark about how
our method is also applicable
to study a coupling by parallel transport.

\subsection{Existing arguments on coupling methods}

As stated above,
we compare our method of the proof
with existing arguments in coupling methods
from a technical point of view.
We hope that the following observation will be helpful
to extend coupling arguments other than our own
in this time-inhomogeneous case.

In order to go into details,
let us review a heuristic (and common) idea of
the construction of a coupling by reflection
as well as that of the derivation of~(\ref{eqdom}).
Given a Brownian particle $X_1$,
we will construct $X_2$ by determining
its infinitesimal motion
$d X_2 (t) \in T_{X_2 (t)} M$
by using $d X_1 (t) \in T_{X_1 (t)} M$.
First we take a minimal $g(t)$-geodesic
$\gm$ joining $X_1 (t)$ and $X_2 (t)$.
Next, by using the parallel transport along $\gm$
associated with the $g(t)$-Levi--Civita connection,
we bring $d X_1 (t)$ into $T_{X_2 (t)} M$.
Finally we define $d X_2 (t)$
as a reflection of it with respect to a hyperplane
being $g(t)$-perpendicular to $\dot{\gm}$
in $T_{X_2 (t)} M$.
From this construction,
the It\^{o} formula implies that
$d_{g(t)} ( X_1 (t) , X_2 (t) )$
should become a semimartingale at least until
$( X_1 (t) , X_2 (t) )$ hits the $g(t)$-cutlocus $\Cut_{g(t)}$.
The semimartingale decomposition is
given by variational formulas of arc length.
On the bounded variation part,
there appear
the time-derivative of $d_{g(t)}$
and (a trace of) the second variation of~$d_{g(t)}$,
which is dominated in terms of the Ricci curvature.
With the aid of our condition~(\ref{eqc-b}),
these two terms are compensated and
a nice domination of
the bounded variation part follows.
Thus the hitting time to 0 of $d_{g(t)} ( X_1 (t) , X_2 (t) )$,
which is the same as $\tau^*$, can be
estimated by that of the dominating semimartingale.
Indeed, we can regard $2 B (t) + d_{g(T_1)} (x_1 , x_2 )$
which appeared in the right-hand side of~(\ref{eqdom})
as the dominating semimartingale.
The effect of our reflection appears
in the martingale part $2 B (t)$
which makes it possible
for the dominating semimartingale to hit $0$.
This construction seems to work
as long as $( X_1 (t) , X_2 (t) )$ is not in the cutlocus.
Moreover, if we succeed
in constructing it beyond the cutlocus,
then the same domination should hold.
Indeed,
the effect of singularity at the cutlocus
should decrease $d_{g(t)} ( X_1 (t) , X_2 (t) )$.
Thus a ``local time at the cutlocus''
will be nonpositive, and hence negligible.

After this observation,
we can conclude that
almost all technical difficulties are
concentrated on the treatment of singularity at the cutlocus
in order to make the heuristic argument rigorous.
In fact,
Theorem~\ref{thmain} is shown in~\cite{Philisem}
by using SDE methods
under the assumption that
the $g(t)$-cutlocus is empty
for every $t \in[ T_1 , T_2 ]$.
It should be remarked that
the joint distribution of the coupled particle
$( X_1 (t) , X_2 (t) )$ could be singular
to the Riemannian measure on $M \times M$
(at least it is the case when $M$ is a flat Euclidean space).
Thus it is not clear that the cutlocus is really ``small''\vadjust{\goodbreak}
for the coupled particle
despite the fact that the cutlocus (as a subset of $M \times M$)
has null $g(t)$-Riemannian measure.

In our approach,
we first construct a coupling of geodesic random walks
and then take a limit to obtain the desired coupling.
Since we first derive a dominating semimartingale
for coupled geodesic random walks,
we need only a \textit{difference inequality}
instead of the It\^o formula.
By virtue of this difference,
we can obtain a desired estimate beyond the cutlocus
by dividing a minimal geodesic joining particles
into small pieces so that the endpoints of
each piece are uniformly away from the cutlocus; see Lemma~\ref{lemc-2var}.
As a~result,
we can avoid extracting a local time at the cutlocus
and directly obtain a dominating process
which does not involve such a term.
Moreover, the dependency on time parameter of the cutlocus
does not cause much difficulty in our approach.

In the time-homogeneous case, there are several arguments
\cite{Crans,Hsu,Wang94,Wang97,Wangbook05}
to construct a coupling by reflection
by approximating it with ones which move
as mentioned above, if they are distant from the cutlocus
and move independently if they are close to the cutlocus.
In some of those arguments,
we need to estimate the size of the total time
when particles are close to the cutlocus.
In such a case, an extension of these arguments
to the time-inhomogeneous case
does not seem straightforward
since the $g(t)$-cutlocus depends on time and
estimates should be more complicated.
The argument in~\cite{Wangbook05} uses supermartingales
to extract the local time at the cutlocus in an implicit way, and
no estimate of times spent around the cutlocus is necessary.
Thus it seems possible to extend his argument
in the time-inhomogeneous case.
Since his argument relies on some detailed properties of parabolic PDEs,
we need to develop time-inhomogeneous analogs of them
to complete this plan.
The fact that our assumption~(\ref{eqc-b}) [or~(\ref{eqc-b0})]
does \textit{not} imply any time-uniform lower bound
of the Ricci curvature by a constant might be an obstacle.

If we employ the theory of optimal transportation,
we will work on couplings of heat distributions
instead of coupling of Brownian motions.
Once we move to the world of heat distributions,
we can expect that the cutlocus is treated more easily
since they are of measure zero
with respect to the Riemannian measure.
However, at this moment,
the theory of optimal transport is not
so strong a tool in this context
for the following two reasons.
First, the range of the theory is restrictive
in the sense that it only deals with couplings
corresponding to the coupling by parallel transport.
Second, the theory of optimal transportation
provides a weaker result than a probabilistic approach does,
even in studying couplings by parallel transport; for instance, see
\cite{McC-ToppWass-RF}
and compare it with~\cite{Arn-Coul-Thalhoriz}.
It should be remarked that
such a difference
between these two approaches
exists even in the time-homogeneous case.

Arnaudon, Coulibaly and Thalmaier
\cite{Arn-Coul-Thalhoriz} recently
developed a new method to construct a coupling,
which works\vadjust{\goodbreak} even in the time-inhomogeneous case.
They consider a one-parameter family of
coupled particles along a curve.
Intuitively speaking,
they concatenate coupled particles along a curve
by iteration of making a coupling by parallel transport.
Since ``adjacent'' particles are
infinitesimally close to each other,
we can ignore singularities on the cutlocus
when we construct a coupled particle
from an ``adjacent'' one.
It should be noted that
their method does not seem to be able to
be applied directly
in order to construct a coupling by reflection.
Indeed, their construction of a chain of coupled particles
heavily relies on a multiplicative (or semigroup) property
of the parallel transport.
However, our reflection operation obviously
fails to possess such a multiplicative property.
Since our reflection map changes orientation,
there is no chance to interpolate it
with a continuous family of isometries.

%
\section{Properties on time-dependent metric}
\label{sec2}

As in Section~\ref{sec1},
let $M$ be a~$m$-dimensional manifold and
$( g(t) )_{t \in[ T_1 , T_2 ]}$
a family of complete Riemannian metrics on $M$
which smoothly depends on~$t$,
for $-\infty< T_1 < T_2 < \infty$.
\begin{remark}
\label{remtime-parameter}
It seems to be restrictive that
our time parameter only runs over
the compact interval $[ T_1 , T_2 ]$.
An example of $g(t)$ we have in mind is
a solution to the backward Ricci flow equation.
In this case, we can work
on a semi-infinite interval $[ T_1 ,\infty)$
only when we study an ancient solution of the Ricci flow.
Thus $T_2 < \infty$ is not so restrictive.
In addition, we could extend
our results to the case on $[ T_1 , \infty)$
with a small modification of our arguments.
It would be helpful to study an ancient solution.
To deal with a singularity of Ricci flow,
it could be nice to work
on a semi-open interval $( T_1 , T_2 ]$,
where $T_1$ is the first time
when a singularity emerges.
In that case, we should be more careful
since we cannot give
``an initial condition at $T_1$''
to define a $g(t)$-Brownian motion on $M$.
\end{remark}

We collect some notation
which will be used in the sequel.
Throughout this paper,
we fix a reference point $o \in M$.
Let $\N_0$ be nonnegative integers.
For $a, b \in\R$,
$a \wg b$ and $a \vee b$ stand
for $\min\{ a, b \}$ and $\max\{ a , b \}$,
respectively.
Let $\Cut_{g(t)} (x)$ be
the set of the $g(t)$-cutlocus of $x$ on $M$.
Similarly, the $g(t)$-cutlocus $\Cut_{g(t)}$ and
the space--time cutlocus $\Cutst$ are
defined by
\begin{eqnarray*}
\Cut_{g(t)}
& : = &
\bigl\{ (x,y) \in M \times M | y \in\Cut_{g(t)} (x) \bigr\},
\\
\Cutst
& : = &
\bigl\{ (t,x,y) \in[ T_1 , T_2 ] \times M \times M | (x,y) \in
\Cut_{g(t)} \bigr\}.
\end{eqnarray*}
Set $D(M) := \{ (x , x) | x \in M \}$.
The distance function with respect to $g(t)$
is denoted by $d_{g(t)} (x,y)$.
Note that
$\Cutst$ is closed
and that $d_{g(\cdot)} ( \cdot, \cdot)$ is
smooth on
$
[ T_1 , T_2 ] \times M \times M
\setminus
( {\Cutst}\cup{[ T_1 , T_2 ]} \times D(M) )
$; see~\cite{McC-ToppWass-RF}; cf.~\cite{K-Phili}.
We denote an open \mbox{$g(s)$-ball} of radius $R$
centered at $x \in M$
by $B^{(s)}_R (x)$.
Some additional notation will be given
at the beginning of the next section.

In the following three lemmas
(Lemmas~\ref{lemmetriccontrol}--\ref{lembddcpt}),
we discuss a local comparison between
$d_{g(t)}$ and $d_{g(s)}$ for $s \neq t$.
Those will be a geometric basis of
further arguments.
%
%
\begin{lem} \label{lemmetriccontrol}
Let $M_0$ be a compact subset of $M$.
Then there exists $\kp= \kp( M_0 )$
such that
\[
\mathrm{e}^{ - 2 \kp| t-s | } g (s)
\le
g (t)
\le
\mathrm{e}^{ 2 \kp| t-s | } g (s)
\]
holds on $M_0$ for $t,s \in[ T_1 , T_2 ]$.
In particular,
if
a minimal $g(s)$-geodesic $\gm$ joining $x,y \in M_0$
is included in $M_0$,
then, for $t \in[ T_1 , T_2 ]$,
\[
d_{g(t)} (x,y) \le\mathrm{e}^{\kp| t-s |} d_{g(s)} (x,y).
\]
\end{lem}
\begin{pf}
Let $\pi \dvtx TM \to M$ be a canonical projection.
Let us define $\hat{M}_0$ by
\[
\hat{M}_0
: =
\bigl\{ ( t, v ) \in[ T_1 , T_2 ] \times TM | \pi(v) \in
M_0 , | v |_{g(t)} \le1 \bigr\}.
\]
Note that $\hat{M}_0$ is closed
since $g(\cdot)$ is continuous.
We claim that $\hat{M}_0$ is sequentially compact.
Let us take a sequence
$( ( t_n , v_n ) )_{n \in\N} \subset\hat{M}_0$.
We may assume $t_n \to t \in[ T_1 , T_2 ]$
and $\pi(v_n) \to p \in M_0$ as $n \to\infty$
by taking a subsequence if necessary.
Let $U$ be a neighborhood of $p$ such that
$\{ v \in TM | \pi(v ) \in U \}
\simeq
U \times\R^m$.
For sufficiently large $n$,
$v_n$ is in $U \times\R^m$ and
we write $v_n = ( p_n , \tilde{v}_n )$.
If we cannot take any convergent subsequence of
$( v_n )_{n \in\N}$,
then $| \tilde{v}_n | \to\infty$
as $n \to\infty$,
where \mbox{$| \cdot|$} stands for
the standard Euclidean norm on $\R^m$
[irrelevant to $( g(t) )_{t \in[ T_1 , T_2 ]}$].
Set $v_n' = ( p_n , | \tilde{v}_n |^{-1} \tilde{v}_n )$.
Then, there exists a subsequence
$( v_{n_k}' )_{k \in\N} \subset( v_n' )_{n \in\N}$
such that
$v_{n_k}' \to v_{\infty}' = (p, \bar{v}' )$
as $n \to\infty$
for some $\bar{v}' \in\R^m$ with $| \bar{v}' | =1$.
Since $g(\cdot)$ is continuous,
$
g( t_{n_k} ) ( v_{n_k}' , v_{n_k}' )
\to
g(t)( v_{\infty}' , v_{\infty}' )
$
as $k \to\infty$.
On the other hand,
$
g ( t_{n_k} ) ( v_{n_k}' , v_{n_k}' )
\le
| \tilde{v}_{n_k} |^{-2} \to0
$
since $g (t_n ) ( v_n , v_n ) \le1$.
Thus $\bar{v}'$ must be 0.
It contradicts with $| \bar{v}'| = 1$.
Hence $\hat{M}_0$ is sequentially compact.

Since
$
\hat{M}_0 \ni( t, v )
\mapsto
\partial_t g (t) (v,v)
$
is continuous,
there exists a constant $\kp= \kp( M_0 ) >0$
such that
$
| \partial_t g (t) ( v, v ) |
\le2 \kp
$
for every $( t , v ) \in\hat{M}_0$.
Take $v \in\pi^{-1}( M_0 )$, $v \neq0_{\pi(v)}$.
Then
\[
\partial_t g (t) ( v, v )
= | v |_{g(t)}^2
\partial_t g (t)
\bigl( | v |_{g(t)}^{-1} v , | v |_{g(t)}^{-1} v \bigr)
\le 2 \kp| v |_{g(t)}^2 .
\]
Thus\vspace*{1pt} $\partial_t \log g (t) (v,v) \le2 \kp$ holds.
By integrating it from $s$ to~$t$ with $s < t$,
we obtain
$g(t) (v,v) \le\mathrm{e}^{2 \kp( t-s )} g (s) (v,v)$.
We can obtain the other inequality similarly.

For the latter assertion,
for $a,b$ with $\gm(a) = x$ and $\gm(b) = y$,
\begin{eqnarray*}
d_{g(t)} (x,y)
&\le&
\int_a^b | \dot{\gm}(u) |_{g(t)} \,du
\le
\mathrm{e}^{\kp| t-s | }
\int_a^b | \dot{\gm}(u) |_{g(s)} \,du
\\
&=&
\mathrm{e}^{\kp|t-s|} d_{g(s)} (x,y).
\end{eqnarray*}
\upqed\end{pf}
\begin{lem} \label{lemballcontrol}
For $R > 0$, $x \in M$ and $t \in[ T_1 , T_2 ]$,
there exists
$\delta= \delta( x, t,\allowbreak R ) > 0$ such that
$\bar{B}^{(s)}_{r} (x) \subset\bar{B}^{(t)}_{3r} (x)$
for $r \le R$ and $s \in[ T_1 , T_2 ]$
with $| s - t | \le\delta$.
\end{lem}

\begin{pf}
Set $\kp:= \kp( \bar{B}^{(t)}_{3R} (x) )$ as in Lemma
\ref{lemmetriccontrol} and $\delta: = \kp^{-1} \log2$. Take $p
\in\bar{B}^{(s)}_{r} (x)$ and a minimal $g(s)$-geodesic $\gm \dvtx
[a,b] \to M$ joining $x$ and $p$. Suppose that there exists $u_0
\in[a,b]$ such that $\gm(u_0) \in\bar{B}^{(t)}_{3r} (x)^c$. Let $
\bar{u}_0 : = \inf \{ u \in[ a,b ] | \gm( u ) \in\bar{B}^{(t)}_{3r}
(x)^c \} $. Since $ \gm( [ a , \bar{u}_0 ] ) \subset \bar{B}^{(t)}_{3r}
(x) \subset \bar{B}^{(t)}_{3R} (x) $ and $d_{g(t)} ( x ,\break \gm( \bar{u}_0
) ) = 3r$, Lem\-ma~\ref{lemmetriccontrol} yields
\[
d_{g(s)} (x,p)
\ge
\int_a^{\bar{u}_0} | \dot{\gm}(u) |_{g(s)} \,du
\ge
\mathrm{e}^{-\kp\delta}
\int_a^{\bar{u}_0} | \dot{\gm}(u) |_{g(t)} \,du
\ge
\frac{3r}{2}.
\]
This is absurd.
Hence $\gm([a,b]) \in\bar{B}^{(t)}_{3r} (x)$.
In particular,
$\gm(b) = p \in\bar{B}^{(t)}_{3r} (x)$.
\end{pf}
\begin{lem} \label{lembddcpt}
For $R > 0$,
there exists a compact subset
$M_0 = M_0 (R)$ of~$M$
such that
%
\begin{equation} \label{eqbddcpt}
\Bigl\{ p \in M \big| \inf_{t \in[ T_1 , T_2 ] } d_{g(t)} (
o , p ) \le R \Bigr\}
\subset M_0 .
\end{equation}
\end{lem}
\begin{pf}
For each $t \in[ T_1 , T_2 ]$,
take $\delta( o, t, R + 1 ) > 0$
according to Lem\-ma~\ref{lemballcontrol}.
Take
$\{ t_i \}_{i=1}^n \subset[ T_1 , T_2 ]$
such that
\[
[ T_1 , T_2 ] \subset \bigcup_{i=1}^n \bigl( t_i - \delta( o, t_i , R + 1) ,
t_i + \delta( o, t_i , R + 1) \bigr) .
\]
Let us define a compact set $M_0 \subset M$ by $ M_0 := \bigcup_{i=1}^n
\bar{B}^{(t_i)}_{3R} (o) $. Take \mbox{$p \in M$} such that $\inf_{T_1 \le t
\le T_2} d_{g(t)} (o,p) \le R$. For $\ep\in(0,1)$, take $s \in[ T_1 ,
T_2 ]$ such that $d_{g(s)} (o,p) \le R + \ep$. Then there exists $j
\in\{ 1, \ldots, N \}$ such that $|s - t_j | < \delta( o, t_j , R +
1)$. By Lemma~\ref{lemballcontrol}, it implies $ p \in\bar{B}^{(s)}_{R
+ \ep} (o) \subset \bar{B}^{(t_j)}_{3( R + \ep)} (o) \subset\break
\bigcup_{i=1}^n \bar{B}^{(t_i)}_{3( R + \ep)} (o) $. Hence the
conclusion follows by letting $\ep\downarrow0$.
\end{pf}

Another useful consequence of
Lemmas~\ref{lemmetriccontrol} and~\ref{lemballcontrol}
is the following:
\begin{lem} \label{lemd-conti}
$d_{g(\cdot)} (\cdot, \cdot)$ is continuous
on $[T_1 , T_2 ] \times M \times M$.
\end{lem}
\begin{pf}
Since the topology
on $[T_1 ,T_2 ] \times M \times M$
is metrizable,
it suffices to show
$
\lim_{n\to\infty} d_{g(t_n)} ( x_n , y_n )
=
d_{g(t)} (x,y)
$ when $(t_n , x_n , y_n ) \to( t , x , y )$
as $n \to\infty$.
By the triangle inequality,
%
\begin{eqnarray} \label{eqd-conv}
\bigl| d_{g(t_n)} ( x_n , y_n ) - d_{g(t)} ( x , y )
\bigr|
&\le&
\bigl| d_{g(t_n)} (x,y) - d_{g(t)} (x,y) \bigr|
\nonumber\\[-8pt]\\[-8pt]
&&{} + d_{g(t_n)} (x,x_n)+ d_{g(t_n)} (y,y_n).
\nonumber
\end{eqnarray}
Take $R > 0$ so that $B_{R}^{(t)} (x)$ includes
a minimal $g(t)$-geodesic joining $x$ and $y$.
Take $\kp= \kp( \bar{B}_{4R}^{(t)} (x) )$
according to Lemma~\ref{lemmetriccontrol}.
We can easily see that
every minimal $g(t)$-geodesic
joining $y$ and $y_n$ is included in $B_{2R}^{(t)} (x)$
for sufficiently large $n \in\N$.
Thus Lemma~\ref{lemmetriccontrol} yields
\[
\limsup_{n\to\infty} d_{g(t_n)} (y,y_n)
\le
\limsup_{n \to\infty}
\mathrm{e}^{\kp| t - t_n|} d_{g(t)} (y,y_n)
= 0.
\]
We can show $d_{g(t_n)} (x,x_n) \to0$ similarly.
Take a minimal $g(t_n)$-geodesic
$\gm_n \dvtx [ a$, $b ] \to M$
joining $x$ and $y$.
By our choice of $R$,
Lemma~\ref{lemmetriccontrol} again yields
\[
d_{g(t_n)} ( x, \gm_n (u) )
\le
d_{g(t_n)} ( x, y )
\le
\mathrm{e}^{\kp| t - t_n | } d_{g(t)} (x,y)
\le
\mathrm{e}^{\kp| t - t_n | } R.
\]
It implies
$\limsup_{n\to\infty} d_{g(t_n)} (x,y) \le d_{g(t)} (x,y)$.
In addition,
$\gm_n$ is included in $B_{4R/3}^{(t_n)} (x)$
for sufficiently large $n$.
Thus Lemmas~\ref{lemballcontrol} and~\ref{lemmetriccontrol}
yield
$d_{g(t)} (x,y) \le\mathrm{e}^{\kp| t - t_n |} d_{g(t_n)} (x,y)$.
Hence the conclusion follows
by combining these estimates with~(\ref{eqd-conv}).
\end{pf}

Before closing this section,
we will provide a local lower bound
of the injectivity radius
which is uniform in time parameter.
\begin{lem} \label{leminject}
For every $M_1\! \subset\! M$ compact,
there is $\tilde{r}_0 \!= \!\tilde{r}_0 ( M_1 ) > 0$
such that
$d_{g(t)} (y,z)\! < \!\tilde{r}_0$ implies
$(t,y,z) \!\notin\!\Cutst$
for any
\mbox{$(t,y,z)\! \in\![ T_1 , T_2 ]\! \times\! M_1 \!\times\! M_1$}.
\end{lem}
\begin{pf}
Take $R > 1$ so that $ \sup_{ t \in[ T_1 , T_2 ] } \sup_{ x \in M_1 }
d_{g(t)} ( o , x ) < R - 1 $. By Lemma~\ref{lembddcpt}, there exists a
compact set $M_0 \subset M$ such that~(\ref{eqbddcpt}) holds. For every
$t \in[ T_1 , T_2 ]$ and $x \in M_1$, $(t,x,x) \notin\Cutst$. It
implies that there is $\h_{t,x} \in( 0, 1 )$ such that $( s, y, z )
\notin\Cutst$ whenever\vspace*{1pt} $ d_{g(t)} (x,y) \vee d_{g(t)} (x,z) \vee |t-s|
< \h_{t,x} $ since $\Cutst$ is closed. Thus there exist $N \in\N$ and
$( t_i , x_i ) \in[ T_1 , T_2 ] \times M_1$ ($i = 1, \ldots, N$) such
that
\[
[ T_1 , T_2 ] \times M_1 \subset \bigcup_{i=1}^N \biggl( t_i -
\frac{\h_{t_i , x_i}}{2} , t_i + \frac{\h_{t_i , x_i}}{2} \biggr)
\times B_{\h_{t_i , x_i} / 2 }^{(t_i)} (x_i) .
\]
Set $\tilde{r}_0 > 0$ by
\[
\tilde{r}_0 := \frac12 \exp \biggl( - \frac{\kp}{2} \max_{1 \le i \le N}
\h_{t_i , x_i} \biggr) \min_{1 \le i \le N} \h_{t_i , x_i} ,
\]
where $\kp= \kp( M_0 ) > 0$ is as in Lemma~\ref{lemmetriccontrol}. Take
$(s,y,z) \in[ T_1 , T_2 ] \times M_1 \times M_1$ with $d_{g(s)} ( y, z
) < \tilde{r}_0$. Take $j \in\{ 1 , \ldots, N \}$ so that $ | s - t_j |
\vee d_{g(t_j)} ( x_j , y ) < \h_{t_j , x_j} / 2 $. By virtue of the
choice of $R$ and $M_0$, Lemma~\ref{lembddcpt} yields that every
$g(s)$-geodesic joining $y$ and $z$ is included in $M_0$. Thus Lemma
\ref{lemmetriccontrol} yields
\[
d_{g(t_j)} ( y , z ) \le \mathrm{e}^{\kp| s - t_j |} d_{g(s)} ( y , z )
< \frac{\h_{t_j , x_j}}{2}.
\]
It implies $ |s - t_j| \vee d_{g(t_j)} ( x_j , y ) \vee d_{g(t_j)} (
x_j , z ) < \h_{t_j , x_j} $ and hence $(s,y,z) \notin\Cutst$.
\end{pf}

\section{Approximation via geodesic random walks}
\label{sec3}

Let $( Z(t) )_{t \in[ T_1 , T_2 ]}$ be a~family of smooth vector fields
continuously depending on the parameter $t \in[ T_1 , T_2 ]$. Let
$X(t)$ be the diffusion process associated with the time-dependent
generator $\sL_t = \Dl_{g(t)}/2 + Z(t)$; see~\cite{CoulgtBM} for a
construction of $X(t)$ by solving a SDE on the frame bundle. Note that
$( t , X(t) )$ is a unique solution to the martingale problem
associated with $\partial_t + \sL_\cdot$ on $[ T_1 ,T_2 ] \times M$;
see~\cite{Hsu} for the time-homogeneous case. Its extension to
time-inhomogeneous case is straightforward; see~\cite{Str-Var} also.

In what follows,
we will use several notions in Riemannian geometry
such as exponential map $\exp$, Levi--Civita connection $\nab$,
Ricci curvature $\Ric$ etc.
To clarify the dependency on the metric $g(t)$,
we put $(t)$ on superscript or~$g(t)$ on subscript.
For instance, we use the following symbols:
$\exp^{(t)}$, $\nab^{(t)}$ and $\Ric_{g(t)}$.
We refer to~\cite{Chavel2}
for basics in Riemannian geometry
which will be used in this paper.

For each $t \in[ T_1 , T_2 ]$,
we fix a measurable section
$\Ph^{(t)} \dvtx M \to\mathscr{O}^{(t)} (M)$
of the $g(t)$-orthonormal frame bundle
$\mathscr{O}^{(t)} (M)$ of $M$.
Take a sequence of independent, identically distributed
random variables $\{ \xi_n\}_{n\in\N}$ which are
uniformly distributed on the unit disk in $\R^m$.
Given $x_0 \in M$,
let us define
a~continuously-interpolated
geodesic random walk
$( X^\alpha( t ) )_{t \in[ T_1 , T_2 ]}$
on~$M$
starting from $x_0$
with a scale parameter $\alpha> 0$
inductively.
Let $t_n^{(\alpha)} : = ( T_1 + \alpha^2 n ) \wg T_2$
for $n \in\N_0$.
For $t = T_1 = t_0^{(\alpha)}$,
set $X^\alpha(T_1) := x_0$.
After $X^\alpha(t)$ is defined
for $t \in[ T_1 , t_n^{(\alpha)} ]$,
we extend it to $t \in[ t_n^{(\alpha)} , t_{n+1}^{(\alpha)} ]$
by
\begin{eqnarray*}
\tilde{\xi}_{n+1} & : = & \sqrt{m+2} \Ph^{( t_n^{(\alpha)} )} \bigl(
X^\alpha\bigl( t_n^{(\alpha)} \bigr) \bigr) \xi_{n+1} ,
\\ 
X^\alpha( t ) & := & \exp_{X^\alpha( t_n^{(\alpha)} )}^{(
t_n^{(\alpha)} )} \biggl( \frac{ t - t_n^{(\alpha)} }{\alpha^2} \bigl(
\alpha\tilde{\xi}_{n+1} + \alpha^2 Z\bigl( t_n^{(\alpha)} \bigr) \bigr) \biggr) .
\end{eqnarray*}
%
%
For later use, we define $ N^{(\alpha)} : = \inf \{ n \in\N_0 |
t_{n+1}^{(\alpha)} - t_n^{(\alpha)} < \alpha^2 \} $. This is the total
number of discrete steps of our geodesic random walks with scale
parameter~$\alpha$.
Set $\sC: = C ([ T_1 , T_2 ] \to M)$
and $\sD: = D ( [ T_1 , T_2 ] \to M)$,
the space of right continuous paths on $M$
parametrized with $[T_1 , T_2 ]$
possessing a left limit at every point.
By using a distance $d_{g(T_1)}$ on $M$,
we metrize $\sC$ and~$\sD$ as usual
so that $\sC$ and $\sD$ become Polish spaces;
see~\cite{Ethier-Kurtz}
for a~distance function on $\sD$, for example.
Set $\sC_1 := C ( [ T_1 , T_2 ] \to[ 0 , \infty) )$.
%
Let us define a~time-dependent
$(0,2)$-tensor field $( \nab Z (t))^{\flat}$
by
\[
( \nab Z (t))^\flat(X,Y) : = \tfrac12 \bigl( \bigl\langle\nab_X^{(t)} Z (t) , Y
\bigr\rangle_{g(t)} + \bigl\langle\nab _Y^{(t)} Z(t) , X \bigr\rangle_{g(t)} \bigr) .
\]

\begin{assumption} \label{assnon-explosion}
There exists a locally bounded
nonnegative measurable function $b$
on $[ 0 , \infty)$ such that:
\begin{longlist}
\item
For all $t \in[ T_1 , T_2 )$,
\[
2 ( \nabla Z (t) )^\flat+ \partial_t g (t) \le \Ric_{g(t)} + b \bigl(
d_{g(t)} ( o , \cdot) \bigr) g(t) .
\]
\item\hypertarget{non-explosion}
For each $C > 0$,
a one-dimensional diffusion process $y_t$ given by
\[
d y_t = d \beta_t + \frac12 \biggl(
C 
+ \int_0^{y_t} b ( s ) \,ds \biggr) \,dt,
\]
where $\beta_t$ is a standard Brownian motion,
does not explode.
(This is the case if and only if
\[
\int_1^\infty \exp\biggl( - \int_1^y \mathbf{b} (z) \,dz \biggr) \int_1^y \exp\biggl(
\int_1^z \mathbf{b} (\xi) \,d \xi\biggr) \,dz \,dy = \infty,
\]
where $\mathbf{b}(y) := C + \int_0^{y} b(s) \,ds$; see, e.g.,
\cite{Ik-Wat}, Theorem VI.3.2.)
\end{longlist}
\end{assumption}

Note that~(\ref{eqc-b}) is a special case of Assumption \ref
{assnon-explosion}.
Now, we are in position to state the main theorem of this paper.
\begin{theorem} \label{thIP}
Under Assumption~\ref{assnon-explosion},
$X^\alpha$ converges in law to $X$ in $\sC$ as \mbox{$\alpha\to0$}.
\end{theorem}

Most of arguments in this section will be
devoted to show the tightness, that is:
\begin{prop} \label{proptight}
$( X^\alpha)_{\alpha\in(0,1)}$ is tight in $\sC$.
\end{prop}

In fact, as we will see in the following,
Proposition~\ref{proptight} easily implies Theorem~\ref{thIP}.
\begin{pf*}{Proof of Theorem~\ref{thIP}}
By virtue of Proposition~\ref{proptight},
for any subsequence of $( X^\alpha)_{\alpha\in(0,1)}$
there exists a further subsequence $ ( X^{\alpha_k} )_{k \in\N}$
which converges in law in $\sC$ as $k \to\infty$.
Thus it suffices to show that
this limit has the same law as $X$.
Let $( \beta^\alpha(t) )_{t \in[0,\infty)}$ be
a Poisson process of intensity $\alpha^{-2}$
which is independent of $\{ \xi_n \}_{n \in\N}$.
Set
\[
\bar{\beta}^\alpha(t) := \bigl( T_1 + \alpha^2 \beta^\alpha( t - T_1 ) \bigr) \wg
t_{N^{(\alpha)}}^{(\alpha)} .
\]
Then the Poisson subordination $X^{\alpha_k} ( \bar{\beta}^{\alpha_k}
(\cdot) )$ also converges in law in $\sD$ to the same limit; see
\cite{Billi}, for instance. Note that $( \bar{\beta}^{\alpha} (t) ,
X^{\alpha} ( \bar{\beta}^{\alpha} (t) ) )_{t \in[ T_1 , T_2 ]}$ is a
time-inhomogeneous Markov process. The associated semigroup
$P^{(\alpha)}_t$ and its generator $\tilde{\sL}^{(\alpha)}$ are given
by
\begin{eqnarray*}
P^{(\alpha)}_t f & := & \mathrm{e}^{ - ( t - T_1 )\alpha^{-2}} \Biggl(
\sum_{l=1}^{ N^{(\alpha)} } \frac{ ( ( t - T_1 ) \alpha^{-2} )^l }{l!}
\bigl(q^{(\alpha)} \bigr)^l f
\\
&&\hphantom{\mathrm{e}^{ - ( t - T_1 )\alpha^{-2}} \Biggl(}
{} + \sum_{l > N^{(\alpha)} } \frac{ ( ( t - T_1 ) \alpha^{-2} )^l
}{l!} \bigl( q^{(\alpha)} \bigr)^{ N^{(\alpha)} } f \Biggr),
\\
\tilde{\sL}^{(\alpha)} f & := & \alpha^{-2} \bigl( q^{(\alpha)}f - f\bigr) ,
\end{eqnarray*}
where
\[
q^{(\alpha)} f (t,x) := \E\bigl[ f \bigl( t + \alpha^2 , \exp^{(t)}_x \bigl(
\alpha\sqrt{m+2} \Ph^{(t)}(x) \xi_1 + \alpha^2 Z ( t ) \bigr) \bigr)\bigr].
\]
We can easily prove
$\tilde{\sL}^{(\alpha)} f \to( \partial_t + \sL_\cdot) f$
uniformly as $\alpha\to0$
for $f \in C^\infty_0 ( [ T_1 ,\break T_2 ] \times M )$.
Since
$(
\bar{\beta}^{\alpha} (t) ,
X^{\alpha} ( \bar{\beta}^{\alpha} (t) )
)_{t \in[ T_1 , T_2 ]}$
is a solution to the martingale problem
associated with $\tilde{\sL}^{(\alpha)}$,
the limit in law of
$(
\bar{\beta}^{\alpha_k} (\cdot) ,
X^{\alpha} ( \bar{\beta}^{\alpha_k} (\cdot) )
)$ 
solves the martingale problem associated
with $\partial_t + \sL_\cdot$.
By the uniqueness of the martingale problem,
this limit has the same law as
that of $( t , X(t) )_{t \in[ T_1 , T_2 ]}$.
It completes the proof.
\end{pf*}
%
%
\begin{remark}
\label{remtight}
Proposition~\ref{proptight} also asserts that
any subsequential limit in law
is a probability measure on $\sC$.
Since we have not added any cemetery point to $M$
in the definition of $\sC$,
Theorem~\ref{thIP} implies that $X$ cannot explode.
It almost recovers the result in~\cite{K-Phili}.
Our assumption is slightly stronger than that in~\cite{K-Phili}
on the point where we require \hyperlink{non-explosion}{(ii)}
for \textit{all} $C > 0$, not a given constant.
Note that
we will use Assumption~\ref{assnon-explosion}\hyperlink{non-explosion}{(ii)}
only for a~specified constant $2 C_0$ given in Lemma~\ref{lemdrift}.
However, its expression looks complicated, and
it seems to be less interesting
to provide an explicit bound.
\end{remark}

Now we introduce some additional notation
which will be used in the rest of this paper.
For $t \in[ T_1 , T_2 ]$,
we define $\lfloor t \rfloor_\alpha$ by
\[
\lfloor t \rfloor_\alpha := \sup\{ \alpha^2 n + T_1 | n \in\N_0 ,
\alpha^2 n + T_1 < t \}.
\]
%
Set $\sF_n : = \sg( \xi_1 , \ldots, \xi_n )$.
For $R > 1$,
let us define
$\sg_R \dvtx \sC_1 \to[ T_1 , T_2 ] \cup\{ \infty\}$
by
\[
\sg_R (w) : = \inf \{ t \in[ T_1 , T_2 ] | w (t) > R - 1 \},
\]
where $\inf\varnothing= \infty$. We write $\hat{\sg}_R := \sg_R ( d_{g
(\cdot)} ( o , X^\alpha(\cdot) ) ) $ and $ \bar{\sg}_R : = \alpha^{-2}
( \lfloor\hat{\sg}_R \rfloor_\alpha- T_1 ) + 1 $. Note that
$\bar{\sg}_R$ is an $\sF_n$-stopping time. For each $t \in[ T_1 , T_2
]$ and $x ,\break y \in M$ with $x \neq y$, we choose a minimal unit-speed
$g(t)$-geodesic $\gm_{xy}^{(t)} \dvtx [ 0,\break d_{g(t)} (x,y) ] \to M$ from $x$
to $y$. Note that we can choose $\gm_{xy}^{(t)}$ so that $(x,y)
\mapsto\gm_{xy}^{(t)}$ is measurable in an appropriate sense; see, for
example,~\cite{Renespoly}. We use the same symbol $\gm_{xy}^{(t)}$ for
its range $\gm_{xy}^{(t)} ( [ 0, d_{g(t)} (x,y)] )$. 

\subsection{A uniform bound for the escape probability}
\label{secUni-prob}

The goal of this subsection is to show the
following:
\begin{prop} \label{propnon-explosion}
$
\lim_{R \uparrow\infty} \limsup_{\alpha\downarrow0}
\PP[ \hat{\sg}_R \le T_2 ] = 0
$.
\end{prop}

For the proof,
we will establish a discrete analog of
a comparison argument
for the radial process
as discussed in~\cite{K-Phili}.
From now on,
we fix $R > 1$ sufficiently large
so that $d_{g(T_1)} ( o , x_0 ) < R-1$
until the final line of the proof of
Proposition~\ref{propnon-explosion}.
We also fix a compact set
$M_0 \subset M$ satisfying~(\ref{eqbddcpt}).
Set $r_0 := \tilde{r}_0 \wg(1/2)$,
where $\tilde{r}_0 = \tilde{r}_0 ( M_0 )$
is as in Lemma~\ref{leminject}.\vadjust{\goodbreak}

The first step for proving Proposition~\ref{propnon-explosion}
is to show
a difference inequality for the radial process
$d_{g(t)} ( o , X^\alpha(t) )$ (Lemma~\ref{lemvariation}).
It will play the role of
the It\^{o} formula for the radial process
in our discrete setting.
We introduce some notation
to discuss how to avoid the singularity of
$d_{g(\cdot)} ( o , \cdot)$
on $\{ o \} \cup\Cut_{g(\cdot)} (o)$.
For $r > 0$,
let us define a set $A_r' , A_r''$ and $A_r$
as follows:
\begin{eqnarray*}
A_r' & : = & \bigl\{ (t,x,y) \in[ T_1 , T_2 ] \times M_0 \times M_0 |
d_{g(t)} (x,x') + d_{g(t)} (y,y')
+| t- t' | \ge r  \\
&&\hspace*{196.5pt}\mbox{ for any $(t', x', y') \in\Cutst
$}\bigr\},
\\
A_r'' & : = & \bigl\{ (t,x,y) \in[ T_1 , T_2 ] \times M_0 \times M_0 |
d_{g(t)} (x,y) \ge r \bigr\},
\\
A_r & := & A_r' \cap A_r'' .
\end{eqnarray*}
Note that $A_r$ is compact
and that
$d_{g(\cdot)} ( \cdot, \cdot)$ is smooth
on $A_r$.
For $t \in[ T_1 , T_2 ]$ and $p \in M$,
let us define $o_p^{(t)} \in M_0$
by
\[
o_p^{(t)} : =
\cases{ \displaystyle \gm_{o p}^{(t)} \biggl( \frac{r_0}{2} \biggr), &\quad if $( t , o
, p ) \notin A_{r_0}'$,\vspace*{2pt}\cr
o, &\quad otherwise.}
\]
For simplicity\vspace*{-3pt} of notation, we denote $o_{X^\alpha( t_n^{(\alpha)}
)}^{( t_n^{(\alpha)} )}$ by $o_n$. Similarly, we use the symbol $\gm_n$
for $\gm_{o_n X^\alpha( t_n^{(\alpha)} )}^{( t_n^{(\alpha)} )}$
throughout this section. Note\vspace*{1pt} that $( t , o_{p}^{(t)} , p )
\notin\Cutst$ holds. Furthermore, it is uniformly separated from
$\Cutst$ in the following sense:
%
\begin{lem} \label{lemaway0}
There exist $r_1 > 0$ and $\delta_1 > 0$
such that the following holds:
let $t_0 , t \in[ T_1 , T_2 ]$
with $t - t_0 \in[ 0 , \delta_1]$.
Let $p_0 \in B_{R-1}^{(t_0)} (o)$ and
$p \in B_{\delta_1}^{(t_0)} (p_0)$.
Then we have:
\begin{longlist}
\item
\hypertarget{approxerror0}
$
d_{g(t)} ( o , p )
\le
\mathrm{e}^{\kp(t - t_0)}
( d_{g(t_0)} ( o , p_0 ) + d_{g(t_0)} ( p_0 , p ) )
$;
\item
\hypertarget{awaycutlocus0}
$( t , o_{p_0}^{(t_0)} , p ) \in A_{r_1}$
when $p_0 \notin B_{r_0}^{(t_0)} (o)$.
\end{longlist}
Here $\kp= \kp( M_0 ) > 0$ is given
according to Lemma~\ref{lemmetriccontrol}.
\end{lem}

By applying Lemma~\ref{lemaway0} to $X^\alpha$,
we obtain the following:
\begin{cor} \label{coraway}
There exist $\alpha_0 > 0$ and
$h \dvtx [ 0 , \alpha_0 ] \to[ 0, 1 ]$
with $\lim_{\alpha\downarrow0} h (\alpha) = 0$
such that the following holds:
for $\alpha\le\alpha_0$, $n \in\N_0$ and
$s,t \in[ t_n^{(\alpha)}, t_{n+1}^{(\alpha)} ]$,
when $n < \bar{\sg}_R$:
\begin{longlist}
\item
\hypertarget{approxerror}
$
d_{g(t)} ( o , X^\alpha(s) )
\le
\mathrm{e}^{ \kp\alpha^2 }
( d_{g( t_n^{(\alpha)} )} ( o , X^\alpha( t_n^{(\alpha)} ) ) +
h (\alpha) )
$;
\item
$( t , o_n , X^\alpha( s ) ) \in A_{r_1}$
when $X^\alpha( t_n^{(\alpha)} ) \notin B^{(t_n^{(\alpha)})}_{r_0} (o)$.
\end{longlist}
Here $r_1$ is the same as in Lemma~\ref{lemaway0}.
\end{cor}
\begin{pf}
Set
$
\bar{Z}
: =
\sup_{t \in[ T_1 , T_2 ], x \in M_0}
| Z(t) |_{g(t)} (x)
$.
Note that we have
\[
d_{g( t_n^{(\alpha)} )} \bigl( X^\alpha\bigl( t_n^{(\alpha)} \bigr) , X^\alpha(t) \bigr)
\le
\sqrt{m+2} \alpha+ \bar{Z} \alpha^2
\]
by the definition of $X^\alpha$.
Take $\alpha_0 > 0$ so that
$\sqrt{m+2} \alpha_0 + \bar{Z} \alpha_0^2 \le\delta_1$
and
$\alpha^2 \le\delta_1$
hold, where $\delta_1$ is as in Lemma~\ref{lemaway0}.
Then the conclusion follows
by applying Lemma~\ref{lemaway0}
with $t_0 = t_n^{(\alpha)}$,
$p_0 = X^\alpha( t_n^{(\alpha)} )$ and
$p = X^\alpha( s )$.
\end{pf}
%
%
\begin{pf*}{Proof of Lemma~\ref{lemaway0}}
We show that
\hyperlink{approxerror0}{(i)} holds with $\delta_1 = 1$.
By the triangle inequality,
the proof is reduced to showing
the following two inequalities:
%
\begin{eqnarray}
\label{eqd-1}
d_{g(t)} ( o , p_0 ) & \le& \mathrm{e}^{\kp( t - t_0 )} d_{g( t_0 )} (
o , p_0 ) ;
\\
\label{eqd-2}
d_{g(t)} ( p_0 , p ) & \le& \mathrm{e}^{\kp(t - t_0)} d_{g( t_0 )} (
p_0, p ) .
\end{eqnarray}
Our condition~(\ref{eqbddcpt}) yields that
$\gm_{op_0}^{(t_0)}$ is included in $M_0$.
Thus Lemma~\ref{lemmetriccontrol} yields~(\ref{eqd-1}).
When $p \in B_1^{(t_0)} (p_0)$,
we have $\gm_{p_0 p}^{(t_0)} \subset B_R^{(t_0)} (o)$.
Hence
(\ref{eqbddcpt}) and Lemma~\ref{lemmetriccontrol}
yield~(\ref{eqd-2})
in a similar way as~(\ref{eqd-1}).

Let us consider \hyperlink{awaycutlocus0}{(ii)}.
For simplicity of notation,
we denote $o_{p_0}^{(t_0)}$ by $o'$
in this proof.
We assume
that $t - t_0 \in[ 0 , \delta]$ and
$p \in B_\delta^{(t_0)} (p_0)$ hold
for $\delta> 0$.
First we will show
$(t , o' , p ) \in A_{r_0 / 4}''$
when $\delta$ is sufficiently small.
Note that
$( t_0 , o' , p_0 ) \in A_{r_0 / 2}''$ holds
since $p_0 \notin B^{(t_0)}_{r_0} (o)$ and
$d_{g( t_0 )} ( o , o' ) \in\{ r_0 / 2 , 0 \}$.
Let $q \in\gm_{o' p_0}^{(t)}$.
By the triangle inequality,
%
\begin{equation}
\label{eqd-3}
d_{g(t)} ( o , q ) \le d_{g(t)} ( o , o' ) + d_{g(t)} (
o' , p_0 ) .
\end{equation}
Since $r_0 /2 < 1 < R$ holds,
(\ref{eqbddcpt}) yields
$\gm^{(t_0)}_{oo'} \subset M_0$
when $o' \neq o$.
We can easily see that
$\gm_{o'p_0}^{(t_0)} \subset\gm_{o p_0}^{(t_0)} \subset M_0$.
Thus,
by applying Lemma~\ref{lemmetriccontrol} to~(\ref{eqd-3}),
%
\begin{eqnarray}\label{eqd-4}
d_{g(t)} ( o , q )
& \le&
\mathrm{e}^{\kp(t - t_0)}
\bigl( d_{g(t_0)} ( o , o' ) + d_{g(t_0)} ( o' , p_0 ) \bigr)
\nonumber\\[-8pt]\\[-8pt]
& \le & (R-1) \mathrm{e}^{\kp\delta}.\nonumber
\end{eqnarray}
Take $\delta_2 := 1 \wg( \kp^{-1} \log(R / (R-1) ) )$.
Then, for any $\delta\in( 0 , \delta_2 )$,
(\ref{eqd-4}) and~(\ref{eqbddcpt}) imply
$\gm_{o'p_0}^{(t)} \subset M_0$.
Hence the triangle inequality,
Lemma~\ref{lemmetriccontrol} and
(\ref{eqd-2}) yield
%
\begin{eqnarray}\label{eqd-5}
d_{g(t)} ( o' , p )
& \ge&
d_{g(t)} ( o' , p_0 )
-
d_{g(t)} ( p_0 , p )
\nonumber\\
& \ge&
\mathrm{e}^{-\kp(t - t_0)}
d_{g(t_0)} ( o' , p_0 )
-
\mathrm{e}^{\kp(t - t_0)}
d_{g(t_0)} ( p_0 , p )
\\
& \ge&
\frac{\mathrm{e}^{-\kp\delta} r_0}{2}
-
\mathrm{e}^{\kp\delta} \delta,\nonumber
\end{eqnarray}
when $\delta\le\delta_2$. Thus there exists $ \delta_3 = \delta_3 (
\kp, r_0 , R) \in (0, \delta_2 ] $ such that the right-hand side of
(\ref{eqd-5}) is greater than $r_0 / 4$ whenever $\delta\in( 0 ,
\delta_3 )$. Hence $( t , o' , p ) \in A_{r_0 / 4}''$ holds in such a
case.\vadjust{\goodbreak}

Next we will show that there exists $r_1' > 0$ such that $( t , o' , p
) \in A_{r_1'}'$ holds for sufficiently small $\delta$. Once we have
shown it, the conclusion holds with $r_1 = r_1' \wg( r_0 / 4 )$. As we
did in showing $(t , o' , p ) \in A''_{r_0 / 4}$, we begin with
studying the corresponding statement for $( t_0 , o' , p_0 )$. More
precisely, we claim that there exists $r_1'' \in( 0, 1 )$ such that $ (
t_0 , o' , p_0 ) \in A_{r_1''} $. When $o' = o$, $( t_0 , o' , p_0 )
\in A_{r_0}'$ directly follows from the definition of $o' =
o_{p_0}^{(t_0)}$. When $o' \neq o$, set
\begin{eqnarray*}
H &: =&
\bigl\{ (t,x,y) \in[ T_1 , T_2 ] \times M_0 \times M_0 |
r_0 \le d_{g(t)} (o,y) \le R-1 ,\\[-2pt]
&&\hspace*{6.4pt}  d_{g(t)} (o,x) =
r_0 / 2 ,d_{g(t)} (x,y) = d_{g(t)} (o,y) - d_{g(t)} (o,x)\bigr\} .
\end{eqnarray*}
Note that $H$ is compact and that $H \cap\Cutst= \varnothing$ holds
since $(t,x,y) \in H$ implies that $x$ is on a minimal $g(t)$-geodesic
from $y$ to $o$. Since $( t_0 , o' , p_0 ) \in H$ by the definition of
$o'$, it suffices to show that there exists $\tilde{r}_1 > 0$ such that
$H \subset A_{\tilde{r}_1}'$. Indeed, the claim will be shown with
$r_1'' = \tilde{r}_1 \wg r_0$ once we have proved it. Suppose that $H
\subset A_r'$ does not hold for any $r \in( 0 , 1 )$. Then there are
sequences $(t_j , x_j , y_j ) \in H$, $(t_j' , x_j' , y_j' )
\in\Cutst$, $j \in\N$, such\vspace*{1pt} that $ | t_j - t_j' | + d_{g(t_j)} (x_j ,
x_j') + d_{g(t_j)} (y_j , y_j') \to0$ as $j \to\infty$. We may assume
that $( (t_j , x_j , y_j ) )_j$ converges. Since $(t_j , x_j , y_j )
\in H$, $x_j' , y_j' \in M_0$ holds for sufficiently large $j$. Thus we
can take a convergent subsequence of $( (t_j' , x_j' , y_j' ) )_j$.
Since\vspace*{1pt} $\Cutst$ and~$H$ are closed, and $d_{g(\cdot)} ( \cdot, \cdot)$
is continuous, it contradicts with \mbox{$H \cap\Cutst= \varnothing$}.

To complete the proof, we show that there exists $\delta_1 \in( 0 ,
\delta_3 ]$ such that $( t , o' , p ) \in A_{r_1'' / 2}'$ when
$\delta\in( 0 , \delta_1 )$. Suppose that\vspace*{-2pt} there exists $(t' , x' , y' )
\in\Cutst$ such that $ | t - t' | + d_{g(t)} ( o' , x' ) + d_{g(t)} ( p
, y' ) < r_1'' / 2 $. For any $q \in\gm_{py'}^{(t)}$, the triangle
inequality and assertion \hyperlink{approxerror0}{(i)} yield
%
\begin{equation} \label{eqd-6}
d_{g(t)} ( o , q ) \le d_{g(t)} ( o , p ) + d_{g(t)} ( p , y' ) \le
\mathrm{e}^{ \kp\delta} ( R - 1 + \delta) + r_1'' / 2 .
\end{equation}
A similar observation implies $d_{g(t)} ( o , q' ) \le(
\mathrm{e}^{\kp\delta} r_0 + r_1'' ) / 2$ for $q' \in\gm_{o'x'}^{(t)}$.
Thus there is $\delta_4 = \delta_4 ( \kp, R ) \in( 0 , \delta_3 ]$ such
that the right-hand side of~(\ref{eqd-6}) is less than $R$ and $(
\mathrm{e}^{\kp\delta} r_0 + r_1'' ) / 2 \le R$ whenever $\delta\in( 0
, \delta_4 )$. In such a case, $\gm_{py'}^{(t)} \subset M_0$ and
$\gm_{o'x'}^{(t)} \subset M_0$ hold. Since $( t_0 , o' , p_0 ) \in
A_{r_1''}'$, Lemma~\ref{lemmetriccontrol} yields
%
\begin{eqnarray} \label{eqd-7}
&&| t - t' |
+ d_{g(t)} ( o' , x' )
+ d_{g(t)} ( p , y' )
\nonumber\\[-2pt]
&&\qquad \ge
| t_0 - t' | - \delta
+ \mathrm{e}^{-\kp\delta} d_{g(t_0)} ( o' , x' )
+ \mathrm{e}^{-\kp\delta} d_{g(t_0)} ( p , y' )
\\[-2pt]
&&\qquad \ge
\mathrm{e}^{- \kp\delta} r_1''
+ ( 1 - \mathrm{e}^{- \kp\delta} ) | t_0 - t' |
- \delta
- \mathrm{e}^{- \kp\delta} \delta.\nonumber
\end{eqnarray}
Take $\delta_1 = \delta_1 ( \kp, r_1'' ) \in( 0 , \delta_4 ]$
so that the right-hand side of~(\ref{eqd-7})
is greater than $r''_1 / 2$ when $\delta\in( 0 , \delta_1 )$.
Then~(\ref{eqd-7}) is absurd
for any $\delta\in( 0 , \delta_1 )$.
Thus it implies the conclusion.\vspace*{-3pt}
\end{pf*}

We prepare some notation
for the second variation formula for the arc length.
Let $\nabla^{(t)}$ be the $g(t)$-Levi--Civita connection
and $\mathcal{R}^{(t)}$ the $g(t)$-curva\-ture tensor\vadjust{\goodbreak}
associated with $\nabla^{(t)}$.
For a smooth curve $\gm$ and
smooth vector fields $U, V$ along~$\gm$,
the index form $I_\gm^{(t)} ( U , V )$ is given by
\[
I_\gm^{(t)} ( U , V ) : = \int_\gm \bigl( \bigl\langle\nab_{\dot{\gm}}^{(t)} U ,
\nab_{\dot{\gm}}^{(t)} V \bigr\rangle_{g(t)} - \bigl\langle\mathcal{R}^{(t)} ( U,
\dot{\gm} ) \dot {\gm}, V \bigr\rangle_{g(t)} \bigr) \,ds .\vspace*{-2pt}
\]
We write $I_\gm^{(t)} (U,U) =: I_\gm^{(t)} (U)$ for simplicity of
notation. Let $G_{t,x,y} (u)$ be the solution to the following initial
value problem on $[ 0, d (x,y) ]$:
\[
\cases{\displaystyle  G_{t,x,y}'' ( u ) = - \frac{ \Ric_{g(t)} ( \dot{\gm}_{xy}^{(t)}
(u) , \dot{\gm}_{xy}^{(t)} (u) ) }{ m-1 } G_{t,x,y} (u) , \vspace*{1pt}\cr
\displaystyle G_{t,x,y}(0) = 0 ,\qquad G_{t,x,y}' (0) = 1 .}\vspace*{-2pt}
\]
Note that $G_{t,x,y} (u) > 0$ for $u \in( 0, d(x,y) ]$
if $y \notin\Cut_{g(t)} (x)$; see~\cite{K-Phili}, proof of Lemma 9.
For simplicity,\vspace*{-2pt}
we write $G_n : = G_{t_n^{(\alpha)} , o_n , X^\alpha( t_n^{(\alpha)} )}$.
When $X^\alpha( t_n^{(\alpha)} ) \notin B_{r_0}^{( t_n^{(\alpha)} )} (o)$,
we define a vector field $V^\dag$ along $\gm_n$
for each $V \in T_{X^\alpha( t_n^{(\alpha)} )} M$
by
\[
V^\dag(\gm_n (u)) : = \frac{ G_n (u) }{ G_n ( d_{g( t_n^{(\alpha)} )} (
o_n , X^\alpha( t_n^{(\alpha)} ) ) ) } \bigl( \DD_{\gm_n}^{(t_n^{(\alpha)})}
V \bigr)( \gm_n (u) ),\vspace*{-2pt}
\]
where $\DD_{\gm_n}^{(t_n^{(\alpha)})} V$ is the parallel vector field
along $\gm_n$ of $V$ associated with $\nab^{( t_n^{(\alpha)} )}$. Take
$v \in\R^m$. By using these notations, for $n \in\N_0$ with $n <
N^{(\alpha)}$, let us define $\lambda_{n+1}$ and $\Lm_{n+1}$ by
\begin{eqnarray*}
\lambda_{n+1}
& : = &
\langle\tilde{\xi}_{n+1} , \dot{\gm}_n \rangle_{g ( t_n^{(\alpha
)} )}
,
\\[-2pt]
\Lm_{n+1}
& : = &
\partial_t d_{g( t_n^{(\alpha)} )} ( o , o_n )
+
\partial_t d_{g( t_n^{(\alpha)} )} \bigl( o_n , X^\alpha\bigl( t_n^{(\alpha)}
\bigr) \bigr)
\\[-2pt]
&&{} + \bigl\langle Z \bigl( t_n^{(\alpha)} \bigr) , \dot{\gm}_n
\bigr\rangle_{g ( t_n^{(\alpha)} )} \bigl( X^\alpha\bigl(
t_n^{(\alpha)} \bigr) \bigr) + \tfrac{1}{2} I_{\gm_n}^{( t_n^{(\alpha)}
)} ( \tilde{\xi}_{n+1}^\dag ),\vspace*{-2pt}
\end{eqnarray*}
when $ X^\alpha( t_n^{(\alpha)} ) \notin B^{( t_n^{(\alpha)} )}_{r_0}
(o) $, and $\lambda_{n+1} := \sqrt{m+2} \langle\xi_{n+1} , v
\rangle_{\R^m}$ and $\Lm_{n+1} := 0$ otherwise.\vspace*{-2pt}
%
%
\begin{lem} \label{lemvariation}
If $n < \bar{\sg}_R \wg N^{(\alpha)}$, $\alpha< \alpha_0$ is small
enough and $X^\alpha( t_n^{(\alpha)} ) \notin B_{r_0}^{( t_n^{(\alpha)}
)} (o)$, then
\[
d_{g( t_{n+1}^{(\alpha)} )} \bigl( o, X^\alpha\bigl( t_{n+1}^{(\alpha)}
\bigr) \bigr) \le d_{ g ( t_n^{(\alpha)} ) } \bigl( o , X^\alpha\bigl(
t_n^{(\alpha)} \bigr) \bigr) + \alpha\lambda_{n+1} + \alpha^2 \Lm_{n+1}
+ o ( \alpha^2 )\vspace*{-2pt}
\]
almost surely,
where $\alpha_0$ is as in Corollary~\ref{coraway}.
In addition,
$o( \alpha^2 )$ is controlled uniformly.\vspace*{-2pt}
\end{lem}
\begin{pf}
By virtue of Corollary~\ref{coraway},
for sufficiently small $\alpha$,
the Taylor expansion together with
the second variation formula
yields
%
\begin{eqnarray}\label{eq2nd}
&& d_{g( t_{n+1}^{(\alpha)} )} \bigl( o_n , X^\alpha\bigl( t_{n+1}^{(\alpha)} \bigr) \bigr)
\nonumber\\
&&\qquad = d_{ g ( t_n^{(\alpha)} ) } \bigl( o_n , X^\alpha\bigl( t_n^{(\alpha)} \bigr) \bigr) +
\alpha\lambda_{n+1} + \alpha^2 \partial_t d_{g( t_n^{(\alpha)} )} \bigl( o_n
, X^\alpha\bigl( t_n^{(\alpha)} \bigr) \bigr)
\nonumber\\[-8pt]\\[-8pt]
&&\qquad\quad{} + \alpha^2 \bigl\langle Z \bigl( t_n^{(\alpha)} \bigr) , \dot{\gm}_n \bigr\rangle_{g (
t_n^{(\alpha)} )} \bigl( X^\alpha\bigl( t_n^{(\alpha)} \bigr) \bigr) + \frac{\alpha^2}{2}
I_{\gm_n}^{( t_n^{(\alpha)} )} ( J_{\tilde{\xi}_{n+1}} )
\nonumber\\
&&\qquad\quad{}
+ o ( \alpha^2 ),\nonumber
\end{eqnarray}
where $J_{\tilde{\xi}_{n+1}}$ is
a $g( t_n^{(\alpha)} )$-Jacobi field
along $\gm_n$
with a boundary value condition
$J_{\tilde{\xi}_{n+1}} ( o_n ) = 0$ and
$J_{\tilde{\xi}_{n+1}} ( X^\alpha( t_n^{(\alpha)} ) ) = \tilde{\xi
}_{n+1}$.
Note that $o ( \alpha^2 )$ can be chosen uniformly
since this expansion can be done
on the compact set $A_{r_1}$, and
every geodesic variation is included in $M_0$.
By the index lemma,\vspace*{-1pt}
we have
$
I_{\gm_n}^{( t_n^{(\alpha)} )} ( J_{\tilde{\xi}_{n+1}} )
\le
I_{\gm_n}^{( t_n^{(\alpha)} )} ( \tilde{\xi}_{n+1}^\dag)
$.
Hence the desired inequality follows when \mbox{$o_n = o$}.
In the case $o_n \neq o$, we have
\begin{eqnarray*}
d_{g( t_{n+1}^{(\alpha)} )} \bigl( o , X^\alpha\bigl( t_{n+1}^{(\alpha)} \bigr) \bigr)
& \le&
d_{g( t_{n+1}^{(\alpha)} )} ( o , o_n )
+
d_{g( t_{n+1}^{(\alpha)} )} \bigl( o_n , X^\alpha\bigl( t_{n+1}^{(\alpha)} \bigr) \bigr),
\\
d_{g( t_n^{(\alpha)} )} \bigl( o , X^\alpha\bigl( t_n^{(\alpha)} \bigr) \bigr)
& = &
d_{g ( t_n^{(\alpha)} )} ( o , o_n )
+
d_{g( t_n^{(\alpha)} )} \bigl( o_n , X^\alpha\bigl( t_n^{(\alpha)} \bigr) \bigr).
\end{eqnarray*}
Note that $(t_n^{(\alpha)} , o, o_n )$ is uniformly away from $\Cutst$
because of our choice of $r_0$ and Lemma~\ref{leminject}.
Therefore the conclusion follows
by combining them with~(\ref{eq2nd}).
\end{pf}

Before turning into the next step,
we show the following two complementary lemmas
(Lemmas~\ref{lemLLN} and~\ref{lemdrift})
which provide a nice control of
the second-order term $\Lm_n$
in Lemma~\ref{lemvariation}.
Set
$\bar{\Lm}_n = \E[ \Lm_n | \sF_{n-1} ]$.
%
\begin{lem} \label{lemLLN}
Let $( a_n )_{ n \in\N_0}$ be a uniformly bounded
$\sF_n$-predictable process.
Then
\[
\lim_{\alpha\to0} \alpha^2 \sup \Biggl\{ \Biggl| \sum_{j=n}^{N+1} a_j ( \Lm_j -
\bar{\Lm}_j ) \Biggr| \bigg| n , N \in\N, n \le N \le N^{(\alpha)} \wg \bar{\sg}_R
\Biggr\} = 0
\]
in probability.
\end{lem}
\begin{pf}
Note that the map
$(t,x,y) \mapsto G_{t,x,y} (d(x,y))$
is continuous on $A_{r_1}$.
Since we have $G_{t,x,y} (d(x,y)) > 0$ on $A_{r_1}$,
there exists $K > 0$ such that
$K^{-1} < G_{t,x,y} (d(x,y)) < K$.
This fact together with Corollary~\ref{coraway} yields
$| \Lm_j |$ and $| \bar{\Lm}_j |$
are uniformly bounded
if $j < \bar{\sg}_R$.
Since
$\sum_{j=1}^n a_j ( \Lm_j - \bar{\Lm}_j )$ is
an $\sF_n$-local martingale and
$\bar{\sg}_R$ is $\sF_n$-stopping time,
the Doob inequality yields
%
\begin{equation} \label{eqDoob}
\lim_{\alpha\to0} \alpha^2 \sup_{0 \le N \le N^{(\alpha)}
\wg\bar{\sg}_R} \Biggl| \sum_{j=1}^{N+1} a_j ( \Lm_j - \bar{\Lm}_j ) \Biggr| = 0
\qquad\mbox{in probability.}
\end{equation}
Here we used
the fact $\lim_{\alpha\to0} \alpha^2 N^{(\alpha)} = T_2 - T_1$.
Note that
\begin{eqnarray*}
&&
\bigcup_{ N = 1 }^{N^{(\alpha)} \wg\bar{\sg}_R} \bigcup_{ n = 0 }^N
\Biggl\{ \alpha^2 \Biggl| \sum_{ j = n + 1 }^{N+1} a_j ( \Lm_j - \bar{\Lm}_j ) \Biggr| >
\delta\Biggr\}
\\
&&\qquad
\subset
\bigcup_{ N = 1 }^{N^{(\alpha)} \wg\bar{\sg}_R}
\bigcup_{ n = 1 }^N
\Biggl\{ \alpha^2 \Biggl| \sum_{ j = 1 }^{n} a_j ( \Lm_j - \bar{\Lm
}_j ) \Biggr| > \frac{\delta}{2} \Biggr\}
\cup
\Biggl\{ \alpha^2 \Biggl| \sum_{ j = 1 }^{N+1} a_j ( \Lm_j - \bar {\Lm}_j
) \Biggr| > \frac{\delta}{2} \Biggr\}
\\
&&\qquad
=
\Biggl\{ \alpha^2 \sup_{0 \le N \le N^{(\alpha)} \wg\bar{\sg}_R}
\Biggl| \sum_{j=1}^{N+1} a_j ( \Lm_j - \bar{\Lm}_j ) \Biggr| >
\frac{\delta}{2} \Biggr\}
.
\end{eqnarray*}
Thus the conclusion follows from~(\ref{eqDoob}).
\end{pf}
\begin{lem} \label{lemdrift}
There exists a deterministic
constant $C_0 > 0$ being independent of
$\alpha$ and $R$ such that the following holds:
\[
\bar{\Lm}_{n+1} \le C_0 + \frac12 \int_0^{d_{g( t_n^{(\alpha)} )} (o ,
X^\alpha( t_n^{(\alpha)} ) )} b ( u ) \,du ,
\]
where $b$ is what appeared in Assumption~\ref{assnon-explosion}.
\end{lem}
\begin{pf}
By using $(m+2) \E[ \langle\xi_n , e_i \rangle\langle\xi_n , e_j
\rangle] = \delta_{ij} $, we obtain
\begin{eqnarray*}
\E\bigl[ I_{\gm_n}^{( t_n^{(\alpha)} )} ( \tilde{\xi}_{n+1}^\dag
) \bigr]
& = &
\sum_{j=2}^m
I_{\gm_n}^{( t_n^{(\alpha)} )}
\bigl( \bigl( \Phi^{( t_n^{(\alpha)} )} \bigl( X^\alpha\bigl( t_n^{(\alpha
)} \bigr) \bigr) e_j \bigr)^\dag\bigr)
\\ 
& = & \frac{ (m-1) G_n' ( d ( o_n , X^\alpha( t_n^{(\alpha)} ) ) ) }{
G_n ( d ( o_n , X^\alpha( t_n^{(\alpha)} ) ) ) } .
\end{eqnarray*}
Note that we have
\begin{eqnarray*}
&&
\bigl\langle Z \bigl( t_n^{(\alpha)} \bigr) , \dot{\gm}_n \bigr\rangle_{g (
t_n^{(\alpha)} )}
\bigl( X^\alpha\bigl( t_n^{(\alpha)} \bigr)\bigr)
-
\bigl\langle Z \bigl( t_n^{(\alpha)} \bigr) , \dot{\gm}_n \bigr\rangle_{g (
t_n^{(\alpha)} )} ( o_n )
\\
&&\qquad =
\int_{0}^{d_{g( t_n^{(\alpha)} )} ( o_n , X^\alpha( t_n^{(\alpha)} ))}
\partial_s
\bigl\langle Z \bigl( t_n^{(\alpha)} \bigr) , \dot{\gm}_n \bigr\rangle_{g (
t_n^{(\alpha)} )} ( \gm_n (s) )\bigg|_{s = u}
\,du
\\
&&\qquad =
\int_{0}^{d_{g( t_n^{(\alpha)} )} ( o_n , X^\alpha( t_n^{(\alpha)} ))}
\bigl\langle\nab_{\dot{\gm}_n}^{( t_n^{(\alpha)} )} Z\bigl( t_n^{(\alpha)}
\bigr), \dot{\gm}_n \bigr\rangle_{g ( t_n^{(\alpha)} )} ( \gm_n (u) )
\,du .
\end{eqnarray*}
Recall that, for $ ( t , x , y ) \notin \Cutst $, we have
\[
\partial_t d_{g(t)} ( x, y )
=
\frac12 \int_0^{d_{g(t)} (x,y)}
( \partial_t g(t) )
\bigl( \dot{\gm}^{(t)}_{xy} (u) , \dot{\gm}^{(t)}_{xy}(u) \bigr)
\,du;
\]
cf.~\cite{McC-ToppWass-RF}, Remark 6.
By combining them with Assumption~\ref{assnon-explosion},
%
\begin{eqnarray} \label{eqexpectation}\quad
\bar{\Lm}_{n+1} & = &
\partial_t d_{g ( t_n^{(\alpha)} )} ( o , o_n )
\nonumber\\
&&{} + \frac12 \int_{0}^{d_{g( t_n^{(\alpha)} )} ( o_n , X^\alpha(
t_n^{(\alpha)} ))}
\partial_t g \bigl( t_n^{(\alpha)} \bigr)
( \dot{\gm}_n (u) , \dot{\gm}_n (u) ) \,du
\nonumber\\
&&{} + \bigl\langle Z \bigl( t_n^{(\alpha)} \bigr) , \dot{\gm}_n \bigr\rangle_{g(
t_n^{(\alpha )} )} \bigl( X^\alpha\bigl( t_n^{(\alpha)} \bigr) \bigr) + \frac{ (m-1) G_n' (
d ( o_n , X^\alpha( t_n^{(\alpha)} ) ) ) }{ 2 G_n ( d ( o_n , X^\alpha(
t_n^{(\alpha)} ) ) ) }
\nonumber\\
& \le& \frac12 \int_{0}^{d_{g( t_n^{(\alpha)} )} ( o , X^\alpha(
t_n^{(\alpha)} ))} b ( u ) \,du
\\
&&{} +
\partial_t d_{g( t_n^{(\alpha)} )} ( o , o_n )
+ \bigl\langle Z \bigl( t_n^{(\alpha)} \bigr) , \dot{\gm}_n \bigr\rangle_{g (
t_n^{(\alpha)} )} ( o_n )
\nonumber\\
&&{} + \frac12 \int_{0}^{d_{g( t_n^{(\alpha)} )} ( o_n , X^\alpha(
t_n^{(\alpha)} ))} \Ric_{g ( t_n^{(\alpha)} )} ( \dot{\gm}_n (u) ,
\dot{\gm}_n (u) ) \,du
\nonumber\\
&&{} + \frac{ (m-1) G_n' ( d ( o_n , X^\alpha( t_n^{(\alpha)} ) ) ) }{
2 G_n ( d ( o_n , X^\alpha( t_n^{(\alpha)} ) ) ) }.\nonumber
\end{eqnarray}
Here we used the fact $ b (u) \ge0$
in the case $o_n \neq o$.
Note that
\[
\int_{0}^{r} \Ric_{g ( t_n^{(\alpha)} )} ( \dot{\gm}_n (u) ,
\dot{\gm}_n (u) ) \,du + \frac{ (m-1) G_n' ( r ) }{ G_n ( r ) }
\]
is nonincreasing as a function of $r$.
Indeed, we can easily verify it by taking
a differentiation.
Set
\begin{eqnarray*}
C_1 &: =& \sup_{t \in[ T_1 , T_2 ]} \sup_{ x \in B_{r_0} ^{(t)} (o) } \Bigl( |
Z ( t ) |_{g(t)} (x)
\\
&&\hphantom{\sup_{t \in[ T_1 , T_2 ]} \sup_{ x \in B_{r_0} ^{(t)} (o) } \Bigl(}
{} + \mathop{\sup_{V \in T_x M}}_{| V |_{g(t)} \le1} \bigl( \partial_t g(t) ( V
, V ) + \bigl| \Ric_{g(t)} ( V , V ) \bigr| \bigr) \Bigr).
\end{eqnarray*}
By virtue of Lemma~\ref{lemmetriccontrol}, $C_1 < \infty$ holds. By
applying a usual comparison argument to $G_n' ( r_0 ) / G_n ( r_0 )$,
we obtain
\begin{eqnarray*}
&&
\int_{0}^{d_{g( t_n^{(\alpha)} )} ( o_n , X^\alpha( t_n^{(\alpha)}
))} \Ric_{g ( t_n^{(\alpha)} )} ( \dot{\gm}_n (u) , \dot{\gm}_n (u) )
\,du \\
&&\quad{} + \frac{ (m-1) G_n' ( d ( o_n , X^\alpha( t_n^{(\alpha)} ) ) ) }{
G_n ( d ( o_n , X^\alpha( t_n^{(\alpha)} ) ) ) }
\\
&&\qquad \le C_1 \bigl( r_0 + \coth( C_1 r_0 ) \bigr).
\end{eqnarray*}
Hence the conclusion
with $C_0 = C_1 ( 1 + 3 r_0 / 4 + \coth( C_1 r_0 ) / 2 )$
follows\break from~(\ref{eqexpectation}).\vadjust{\goodbreak}~%
\end{pf}

In the next step, we will introduce
a comparison process to give a control
of the radial process.
Let us define a function $\varphi$
on $(2r_0 , \infty)$ by
\[
\varphi( r ) : = C_0 + \frac12 \int_0^r b (u) \,du,
\]
where $C_0$ is as in Lemma~\ref{lemdrift}.
Let us define another function $\psi$
on $(2 r_0 , \infty)$ so that
$\psi$ is a locally Lipschitz nonincreasing function
satisfying
$\psi( r ) := 2 (r - 2r_0 )^{-1}$
for $r \in( 2 r_0 , 2 r_0 + 1 ]$
and $\psi( r ) : = 0$ for $r \ge2 r_0 + 2$.
Let us define a comparison process $\ro^\alpha(t)$
taking values in $[0,\infty)$
inductively by
\begin{eqnarray}
\ro^\alpha(T_1) & := & d_{g(T_1)} ( o , x_0 ) + 3 r_0 ,
\nonumber\\
\ro^\alpha(t) & := & \ro^\alpha\bigl( t_n^{(\alpha)} \bigr) + \frac{ t -
t_n^{(\alpha)} }{\alpha^2} \bigl( \alpha\lambda_{n+1}
+ \alpha^2 \bigl( \varphi\bigl( \ro^\alpha\bigl( t_n^{(\alpha)} \bigr) \bigr) + \psi\bigl(
\ro^\alpha\bigl( t_n^{(\alpha)} \bigr) \bigr) \bigr) \bigr),\nonumber\\
&&\eqntext{t \in\bigl[ t_n^{(\alpha)} ,
t_{n+1}^{(\alpha)} \bigr].}
\end{eqnarray}
The term $\psi( \ro^\alpha( t_n^{(\alpha)} ) )$ is inserted to avoid a
difficulty coming from the absence of the estimate in Lemma
\ref{lemvariation} on a neighborhood of $o$. By virtue of this extra
term, $\ro^\alpha(t) > 2 r_0$ holds for all $t \in[ T_1 , T_2 ]$ if
$\alpha$ is sufficiently small. Let~$\hat{\sg}_R'$ and~$\bar{\sg}_R'$
be given by $ \hat{\sg}_R' : = \sg_R ( \ro^\alpha) $ and $ \bar{\sg}_R'
: = \alpha^{-2} ( \lfloor\hat{\sg}_R' \rfloor_\alpha- T_1 ) + 1 $. The
following is a modification of an argument in the proof of~\cite{Hsu},
Theorem~3.5.3, into our discrete setting.
%
\begin{lem} \label{lemcomparison1}
For $\delta>0$,
there exist a family of events $( E_\delta^{\alpha} )_\alpha$
with\break $\lim_{\alpha\to0} \PP[ E_\delta^{\alpha} ] = 1$
and a constant $K (\delta) > 0$
with $\lim_{\delta\to0} K (\delta) = 0$
such that, on $E_\delta^{\alpha}$,
\[
d_{ g(t) } ( o, X^\alpha( t ) ) \le \ro^\alpha(t) + K ( \delta)
\]
for $t \in[T_1 , \hat{\sg}_R \wg\hat{\sg}_R' \wg T_2 ]$
and sufficiently small $\alpha$
relative to $\delta$ and $R^{-1}$.
\end{lem}
\begin{pf}
It suffices to show the assertion in the case
$t = t_n^{(\alpha)}$ for some $n \in\N_0$.
Indeed, once we have shown it,
Corollary~\ref{coraway}\hyperlink{approxerror}{(i)} yields
\begin{eqnarray*}
d_{g(t)} ( o , X^\alpha(t) )
& \le& \mathrm{e}^{\kp\alpha^2} \bigl( d_{g(
\lfloor t \rfloor_\alpha)} ( o , X^\alpha( \lfloor t \rfloor_\alpha) )
+ h (\alpha) \bigr)
\\
& \le&
\ro_{ \lfloor t \rfloor_\alpha}^\alpha + K (\delta) + (
\mathrm{e}^{\kp\alpha^2} - 1 ) R + \mathrm{e}^{\kp\alpha^2} h (\alpha)
\\
& \le&
\ro_t^\alpha + K (\delta) + \alpha + ( \mathrm{e}^{\kp\alpha^2}
- 1 ) R + \mathrm{e}^{\kp\alpha^2} h (\alpha)
\end{eqnarray*}
for $t \in[ T_1 , \hat{\sg}_R \wg T_2]$.
Here we used the facts $\varphi\ge0$ and $\psi\ge0$.
From this estimate,
we can easily deduce the conclusion.

For simplicity of notation, we denote $d_{g( t_n^{(\alpha)} )} ( o ,
X^\alpha( t_n^{(\alpha)} ) )$ and $\ro^\alpha( t_n^{(\alpha)} )$
by~$d_n$ and $\ro_n$, respectively, in the rest of this proof. Let us
define a sequence of $\sF_n$-stopping times $S_l$ by $S_0 : = 0$ and
\begin{eqnarray*}
S_{2l+1} & : = & \inf \bigl\{ j \ge S_{2l} | X^\alpha\bigl( t_j^{(\alpha)} \bigr) \in
B_{r_0}^{(t_j^{(\alpha)})} (o) \bigr\} \wg N^{(\alpha)} ,
\\
S_{2l} & : = &
\inf \bigl\{ j \ge S_{2l-1} | X^\alpha\bigl( t_j^{(\alpha)} \bigr)
\notin B_{3r_0 / 2}^{(t_j^{(\alpha)})} (o) \bigr\} \wg N^{(\alpha)} .
\end{eqnarray*}
Since $\ro_n > 2r_0$, it suffices to show the assertion in the case $
S_{2l} \le n < S_{2l+1} \wg\bar{\sg}_R \wg\bar{\sg}_R' $ for some $l
\in\N_0$. Now Lemmas~\ref{lemvariation} and~\ref{lemdrift} imply
\[
d_{j+1} - \ro_{j+1} \le d_{j} - \ro_{j} + \alpha^2 \bigl( \varphi( d_{j} ) -
\varphi( \ro_{j} ) \bigr) + \alpha^2 ( \Lm_{j+1} - \bar{\Lm}_{j+1} ) + o (
\alpha^2 )
\]
for $ j \in [ S_{2l} , S_{2l+1} \wg\sg_R' \wg\bar{\sg}_R' ) $. Here we
used the fact $\psi\ge0$. Let $f_\alpha$ be a~$C^2$-function on $\R$
satisfying:
\begin{longlist}
\item
$f_\alpha|_{( - \infty, -\alpha)} \equiv0$;
$f_\alpha|_{( \alpha, \infty)} (x) = x$;
\item
\hypertarget{convex}
$f_\alpha$ is convex;
\item
\hypertarget{2nderror}
$\alpha^2 \sup_{x \in\R} f_\alpha'' (x) = o (1)$.
\end{longlist}
For example,
a function $f_\alpha$ satisfying these
conditions is constructed by setting
\[
\tilde{f} (x) = \int_{-\infty}^x \int_{-\infty}^t b \exp\biggl( -
\frac{a}{1-s^2} \biggr) 1_{(-1,1)} (s) \,ds \,dt ,
\]
where $a,b$ is chosen to satisfy
\begin{eqnarray*}
\int_{-\infty}^1
\exp\biggl( - \frac{ a }{ 1-s^2 } \biggr) 1_{(-1,1)} (s)
\,ds &=& 1 ,
\\
b \int_{-\infty}^{1}
\int_{-\infty}^t
\exp\biggl( - \frac{ a }{ 1-s^2 } \biggr) 1_{(-1,1)} (s)
\,ds \,dt
&=& 1
\end{eqnarray*}
and $f_\alpha( x ) : = \alpha\tilde{f} ( \alpha^{-1} x )$.
By the Taylor expansion
with condition \hyperlink{2nderror}{(iii)} of $f_\alpha$,
we have
%
\begin{eqnarray} \label{eqTaylor}\quad
f_\alpha( d_{j+1} - \ro_{j+1} ) &\le& f_\alpha( d_{j} - \ro_{j} )
\nonumber\\
&&{} + \alpha^2 f_\alpha' ( d_{j} - \ro_{j} ) \bigl( \varphi( d_{j} ) - \varphi(
\ro_{j} ) + ( \Lm_{j} - \bar{\Lm}_{j} ) \bigr)\\
&&{} + o ( \alpha^2 ) . \nonumber
\end{eqnarray}
Let $C > 0$ be the Lipschitz constant of $\varphi$
on $[ 0, R ]$.
Note that we have
%
\begin{equation} \label{eqTaylor1st}
f_\alpha' ( d_{j} - \ro_{j} ) \bigl( \varphi( d_{j} ) - \varphi( \ro_{j} ) \bigr)
\le C ( d_{j} - \ro_{j} )_+
\end{equation}
since $\varphi$ is nondecreasing.
Now by using~(\ref{eqTaylor}) and~(\ref{eqTaylor1st})
combined with the fact
$d_{S_{2l}} - \ro_{S_{2l}} < - \alpha$
for sufficiently small $\alpha$,
we obtain
%
\begin{eqnarray}\label{eqpre-comparison}\qquad
( d_{n} - \ro_{n} )_+
& \le&
f_\alpha( d_{n} - \ro_{n} )
\nonumber\\
& \le& C \alpha^2 \sum_{j = S_{2k}}^{n-1} ( d_{j} - \ro_{j} )_+
+ \alpha^2 \sum_{j = S_{2k}}^{n-1} f_\alpha' ( d_{j} - \ro_{j} ) (
\Lm_{j+1} - \bar{\Lm}_{j+1} )
\\
&&{} + o (1) .\nonumber
\end{eqnarray}
Here the first inequality follows from condition \hyperlink{convex}{(ii)} of
$f_\alpha$,
and
\mbox{$n \le\alpha^{-2} ( T_2 - T_1 )$} is used
to derive the error term $o(1)$.
Let $E_\delta^{\alpha}$ be an event defined by
\[
E_\delta^{\alpha} : = \Biggl\{ \alpha^2 \sup_{ k \le k' \le N^{(\alpha)}
\wg\bar{\sg }_R } \Biggl| \sum_{j=k}^{k'} f_\alpha' ( d_{j-1} - \ro_{j-1} ) (
\Lm_{j} - \bar{\Lm}_{j} ) \Biggr| < \delta\Biggr\}.
\]
Note that
$a_j = f_\alpha' ( d_{j-1} - \ro_{j-1} )$ is
$\sF_n$-predictable
and uniformly bounded by~$1$.
Thus,
by combining Lemma~\ref{lemLLN} with~(\ref{eqpre-comparison}),
we obtain
\[
( d_{n} - \ro_{n} )_+ \le C \alpha^2 \sum_{j = S_{2l}}^{n-1} ( d_{j} -
\ro_{j} )_+ + 2 \delta
\]
on $E_\delta^{\alpha}$ for sufficiently small $\alpha$.
Thus,
by virtue of a discrete Gronwall inequality
(see~\cite{Will-Wong}, e.g.),
\[
( d_{n} - \ro_{n} )_+ \le 2 \delta \bigl( 1 + ( 1 + C \alpha^2 )^n \bigr) \le 2
\delta\bigl( 1 + \mathrm{e}^{ C (T_2 - T_1)} \bigr).
\]
%
This estimate implies the conclusion.
\end{pf}
\begin{cor} \label{corcomparisonstopping}
For every $R' < R$,
\[
\limsup_{\alpha\to0 } \PP[ \hat{\sg}_R \le T_2 ] \le
\limsup_{\alpha\to0} \PP[ \hat{\sg}_{R'}' \le T_2 ] .
\]
\end{cor}

Now we turn to the proof of our destination
in this section.
%
%
\begin{pf*}{Proof of Proposition~\ref{propnon-explosion}}
By Corollary~\ref{corcomparisonstopping}, the proof of Proposition
\ref{propnon-explosion} is reduced to estimate $\PP[ \hat{\sg}_R' \le
T_2 ]$. To obtain a useful bound of it, we would like to apply the
invariance principle for $\ro^\alpha$. However, there is a~technical
difficulty coming from the unboundedness of the drift term
of~$\ro^\alpha$. To avoid it, we introduce an auxiliary process
$\tilde{\ro}^\alpha$ in the sequel.

Let $\tilde{\varphi}$ be a bounded,
globally Lipschitz
function on $\R$ such that
$\tilde{ \varphi} (r) = \varphi(r) + \psi(r)$
for $r \in[ 2 r_0 + R^{-1} , R ]$.
Let us define an $\R$-valued process
$\tilde{\ro}^\alpha(t)$ inductively by
\begin{eqnarray*}
\tilde{\ro}^\alpha(T_1) & := & d_{g(T_1)} ( o , x_0 ) + 3 r_0 ,
\\
\tilde{\ro}^\alpha(t) & := & \tilde{\ro}^\alpha\bigl( t_n^{(\alpha)} \bigr) +
\frac{ t - t_n^{(\alpha)} }{\alpha^2} \bigl( \alpha\lambda_{n+1} + \alpha^2
\tilde{\varphi} \bigl( \tilde {\ro}^\alpha\bigl( t_n^{(\alpha)} \bigr) \bigr) \bigr),
\qquad t \in\bigl[
t_n^{(\alpha)} , t_{n+1}^{(\alpha)} \bigr] .
\end{eqnarray*}
We also define two diffusion processes $\ro^0 (t)$ and $\tilde{\ro}^0
(r)$ as solutions to the following SDEs:
\begin{eqnarray*}
&& \cases{
d \ro^0 (t) = d B(t) + \bigl( \varphi( \ro^0 (t) ) + \psi( \ro^0 (t) ) \bigr) \,dt
, \vspace*{2pt}\cr \ro^0 (T_1) = d_{g(T_1)} ( o , x_0 ) + 3 r_0 ,}
\\
&&
\cases{ d \tilde{\ro}^0 (t) = d B(t) + \tilde{\varphi} (
\tilde{\ro}^0(t) ) \,d t , \vspace*{2pt}\cr \tilde{\ro}^0 (T_1) = d_{g(T_1)} ( o , x_0
) + 3 r_0 ,}
\end{eqnarray*}
where $( B(t) )_{t \in[ T_1 , T_2 ]}$ is
a standard one-dimensional Brownian motion
with $B(T_1) = 0$.
We claim that $\tilde{\ro}^\alpha$ converges
in law to $\tilde{\ro}^0$ as $\alpha\to0$.
Indeed, we can easily show
the tightness of $( \tilde{\ro}^\alpha)_{\alpha> 0}$
by modifying an argument
for the invariance principle for
i.i.d. sequences since $\tilde{\varphi}$ is bounded.
Then the claim follows
from the same argument
as we used in the proof of Theorem~\ref{thIP}
under Proposition~\ref{proptight}, which is
based on the Poisson subordination and
the uniqueness of the martingale problem.

Let us define $ \h_R \dvtx \sC_1 \to[ T_1 , T_2 ] \cup\{ \infty\} $ by
\[
\h_R (w) : = \inf \{ t \in[ T_1 , T_2 ] | w (t) \le2r_0 + R^{-1} \} .
\]
Then we have
\[
\PP[ \hat{\sg}_R' \le T_2 ] \le \PP [ \sg_R ( \ro^\alpha) \wg\h_R (
\ro^\alpha) \le T_2 ] = \PP [ \sg_R ( \tilde{\ro}^\alpha) \wg\h_R (
\tilde{\ro }^\alpha) \le T_2 ].
\]
Since $\{ w | \sg_R (w) \wg\h_R (w) \le T_2 \}$
is closed in $\sC_1$,
the Portmanteau theorem implies
\begin{eqnarray*}
\limsup_{\alpha\to0} \PP [ \sg_R ( \tilde{\ro}^\alpha) \wg\h_R (
\tilde{\ro }^\alpha) \le T_2 ] &\le& \PP [ \sg_R ( \tilde{\ro}^0 )
\wg\h_R ( \tilde{\ro}^0 ) \le T_2 ]
\\
&=& \PP [ \sg_R ( \ro^0 ) \wg\h_R ( \ro^0 ) \le T_2 ].
\end{eqnarray*}
Since $\ro^0$ is a diffusion process on $( 2 r_0 , \infty)$
which cannot reach the boundary by Assumption~\ref{assnon-explosion},
the conclusion follows.
\end{pf*}

\subsection{Tightness of geodesic random walks}
\label{sectight}

Recall that we have metrized the path space $\sC$
by using $d_{g(T_1)}$.
To deal with
the tightness of $( X^\alpha)_{\alpha\in( 0, 1 )}$
in $\sC$,
we show the following lemma,
which provides a tightness criterion compatible
with the time-dependent metric $d_{g(t)}$.
\begin{lem} \label{lemsufftight}
$( X^\alpha)_{\alpha\in( 0, 1 )}$ is tight
if
%
\begin{eqnarray} \label{eqtight}\qquad
&&\lim_{\delta\to0} \frac{1}{\delta} \limsup_{\alpha\to0} \sup_{
n \in\N_0}
\PP\Bigl[ \sup_{t_n^{(\alpha)} \le s \le( t_n^{(\alpha)} + \delta) \wg T_2}
d_{g(s)} \bigl( X^\alpha\bigl( t_n^{(\alpha)} \bigr) , X^\alpha(s) \bigr) > \ep,
\nonumber\\[-8pt]\\[-8pt]
&&\hspace*{252.6pt}\hat{\sg}_R = \infty \Bigr] = 0 \nonumber
\end{eqnarray}
holds for every $\ep> 0$ and $R > 1$.
\end{lem}
\begin{pf}
$\!\!\!$By following a standard argument
(e.g.,~\cite{Billi}, Theorems~7.3 and~7.4),
we can easily show that
$( X^\alpha)_{\alpha\in(0,1)}$ is tight
if, for every $\ep> 0$,
\[
\lim_{\delta\to0} \frac{1}{\delta} \limsup_{\alpha\to0} \sup_{t \in[
T_1 , T_2 ]} \PP \Bigl[ \sup_{t \le s \le( t + \delta) \wg T_2} d_{g(T_1)} (
X^\alpha(t) , X^\alpha(s) ) > \ep\Bigr] = 0.
\]
Thus, by virtue of Proposition~\ref{propnon-explosion},
$ ( X^\alpha)_{\alpha\in( 0 , 1 )}$ is tight if
\[
\lim_{\delta\to0} \frac{1}{\delta} \limsup_{\alpha\to0} \sup_{t \in[
T_1 , T_2 ]} \PP \Bigl[ \sup_{t \le s \le( t + \delta) \wg T_2} d_{g(T_1)} (
X^\alpha(t) , X^\alpha(s) ) > \ep, \hat{\sg}_R = \infty\Bigr]
= 0\vadjust{\goodbreak}
\]
for every $\ep> 0$ and $R > 1$.
Given $R > 1$,
take $M_0$ and $\kp$
as in Lemmas~\ref{lembddcpt} and~\ref{lemmetriccontrol},
respectively.
Then, for $\ep< 1$ and $s,t \in[ T_1 , T_2 ]$,
\begin{eqnarray*}
&&\bigl\{ d_{g(s)} ( X^\alpha(s) , X^\alpha( \lfloor t \rfloor_\alpha ) )
\le\ep, \hat{\sg}_R = \infty\bigr\}
\\
&&\qquad\subset \bigl\{ d_{g(T_1)} ( X^\alpha(s), X^\alpha( t ) ) \le2 \mathrm
{e}^{\kp( T_2 - T_1 )} \ep, \hat{\sg}_R = \infty\bigr\},
\end{eqnarray*}
if $\alpha$ is sufficiently small.
Thus we have
\begin{eqnarray*}
&&\Bigl\{ \sup_{t \le s \le( t + \delta) \wg T_2} d_{g(T_1)} ( X^\alpha(t) ,
X^\alpha(s) ) > \ep, \hat{\sg}_R = \infty\Bigr\}
\\
&&\qquad\subset \biggl\{ \sup_{\lfloor t \rfloor_\alpha\le s \le( \lfloor t \rfloor
_\alpha+ 2 \delta) \wg T_2} d_{g(s)} ( X^\alpha( \lfloor t \rfloor
_\alpha) , X^\alpha(s) ) > \frac{\mathrm{e}^{-\kp( T_2 - T_1 )} \ep}{2}
, \hat{\sg}_R = \infty\biggr\}
\end{eqnarray*}
for $\alpha^2 \le\delta$, and hence the conclusion follows.
\end{pf}
\begin{pf*}{Proof of Proposition~\ref{proptight}}
Take $R > 1$.
By virtue of Lemma~\ref{lemsufftight},
it suffices to show~(\ref{eqtight}).
Take $M_0 \subset M$ compact
and $\kp$ as in Lemmas~\ref{lembddcpt} and~\ref{lemmetriccontrol}, respectively.
By taking smaller $\ep> 0$,
we may assume that $\ep< \tilde{r_0} / 2$,
where $\tilde{r}_0 = \tilde{r}_0 (M_0)$
is as in Lemma~\ref{leminject}.
Take $n \in\N_0$ with $n < N^{(\alpha)}$.
Let us define a $\sF_k$-stopping time $\zt_\ep$
by
\[
\zt_\ep : = \inf\bigl\{ k \in\N_0 | n \le k \le N^{(\alpha)} , d_{g (
t_k^{(\alpha)} )} \bigl( X^\alpha\bigl( t_n^{(\alpha)} \bigr) , X^\alpha\bigl(
t_k^{(\alpha)} \bigr) \bigr) > \ep\bigr\}.
\]
Then, for sufficiently small $\alpha$,
%
\begin{eqnarray} \label{eqreplace}
&&\Bigl\{ \sup_{ t_n^{(\alpha)} \le s \le( t_n^{(\alpha)} + \delta) \wg T_2 }
d_{g(s)} \bigl( X^\alpha\bigl( t_n^{(\alpha)} \bigr) , X^\alpha( s ) \bigr) \ge2 \ep,
\hat{\sg}_R = \infty\Bigr\}
\nonumber\\[-8pt]\\[-8pt]
&&\qquad\subset \{ \alpha^2 ( \zt_{\ep} - n ) < \delta, \hat{\sg}_R = \infty\}
. \nonumber
\end{eqnarray}
Set $p_k := X^\alpha( t_k^{(\alpha)} )$ for $k \in\N_0$ and $f ( t , x
) := d_{g(t)} ( p_n , x )$. Note that $f^2$ is smooth on $\{ f <
\ep\}$. Let us define $\lambda_k'$ by
\[
\lambda_{k+1}' := \bigl\langle\tilde{\xi}_{k+1} , \dot{\gm}_{ p_n p_k
}^{(t_k^{(\alpha )})} \bigr\rangle_{g (t_k^{(\alpha)})} .
\]
We claim that there exists a constant $C > 0$
such that
%
\begin{equation} \label{eq2ndtight}
f \bigl( t_{k+1}^{(\alpha)} , p_{k+1} \bigr)^2 \le f \bigl( t_k^{(\alpha)} , p_k \bigr)^2 +
2 \alpha f\bigl( t_k^{(\alpha)} , p_k \bigr) \lambda_{k+1}' + C \alpha^2
\end{equation}
for $k \le\zt_\ep\wg N^{(\alpha)}$
on $\{ \hat{\sg}_R = \infty\}$.
Indeed,
in the same way as we did to obtain~(\ref{eq2nd}),
%
\begin{eqnarray} \label{eq2ndtight1}
&&
f \bigl( t_{k+1}^{(\alpha)} , p_{k+1} \bigr)^2\nonumber\\
&&\qquad \le f \bigl( t_k^{(\alpha)} , p_k
\bigr)^2 + 2 \alpha f \bigl( t_k^{(\alpha)} , p_k \bigr) \lambda_{k+1}' + \alpha^2 (
\lambda_{k+1}' )^2
\nonumber\\[-8pt]\\[-8pt]
&&\qquad\quad{}
+ 2 \alpha^2 f\bigl( t_k^{(\alpha)} , p_k \bigr) \bigl(
\partial_t f \bigl( t_k^{(\alpha)} , p_k \bigr)
+ \bigl\langle Z \bigl( t_k^{(\alpha)} \bigr) , \dot{\gm}_{p_n p_k}^{(t_n^{(\alpha
)})} \bigr\rangle_{g ( t_k^{(\alpha)} )} ( p_k ) \bigr)
\nonumber\\
&&\qquad\quad{} + \alpha^2 f\bigl( t_k^{(\alpha)} , p_k \bigr) I_{\gm_{p_n
p_k}^{(t_k^{(\alpha)})}}^{( t_k^{(\alpha)} )} ( J_{\tilde{\xi}_{k+1}} )
+ o ( \alpha^2 ).\nonumber
\end{eqnarray}
Here $o ( \alpha^2 )$ is controlled uniformly.
Let $K_1 > 0$ be a constant satisfying that
the $g(t)$-sectional curvature on $M_0$
is bounded below by $-K_1$ for every $t \in[ T_1 , T_2 ]$.
Such a constant exists since $M_0$ is compact.
Then a comparison argument implies
\[
f\bigl( t_k^{(\alpha)} , p_k \bigr) I_{\gm_{p_n p_k}^{(t_k^{(\alpha)})}}^{(
t_k^{(\alpha)} )} ( J_{\tilde{\xi}_{k+1}} ) \le K_1 f \bigl( t_k^{(\alpha)}
, p_k \bigr) \coth\bigl( K_1 f \bigl( t_k^{(\alpha)} , p_k \bigr) \bigr) .
\]
Here the right-hand side is bounded uniformly
if $k < \zt_\ep\wg N^{(\alpha)}$.
The remaining estimate of the second-order term in
(\ref{eq2ndtight1}) to show~(\ref{eq2ndtight}) is easy
since we are on the event $\{ \hat{\sg}_R = \infty\}$.
Applying~(\ref{eq2ndtight}) repeatedly
from $k=n$ to $k = \zt_\ep$, we obtain
\[
\ep^2 < 2\alpha\sum_{k=n}^{\zt_\ep} f \bigl( t_k^{(\alpha)} , p_k \bigr)
\lambda_{k+1}' + C \delta
\]
on $ \{ \alpha^2 ( \zt_\ep- n ) < \delta, \hat{\sg}_R = \infty\} $. Set
$ N_\delta^{(\alpha)} := \sup\{ k \in\N_0 | k \le\alpha^{-2} \delta+ n
\} $. By taking $\delta< ( 2 C )^{-1} \ep^2$, we obtain
%
\begin{eqnarray}\label{eqreplace2}
&&\{ \alpha^2 ( \zt_\ep- n ) < \delta , \hat{\sg}_R = \infty \}
\nonumber\\
&&\qquad \subset \Biggl\{ \sum_{k=n}^{\zt_\ep} f \bigl( t_k^{(\alpha)} , p_k \bigr) \lambda
_{k+1}' > \frac{\ep^2}{4 \alpha} , \alpha^2 ( \zt_\ep- n ) < \delta,
\hat{\sg}_R = \infty\Biggr\}
\\
&&\qquad \subset \Biggl\{ \sup_{n \le N \le N_\delta^{(\alpha)}} \sum_{k=n}^{N} f \bigl(
t_k^{(\alpha)} , p_k \bigr) 1_{ \{ f ( t_k^{(\alpha)} , p_k ) \le\ep\} }
\lambda_{k+1}' > \frac{\ep^2}{4 \alpha} \Biggr\}.\nonumber
\end{eqnarray}
Set
\[
Y_{k+1} : = \frac{1}{\sqrt{m+2}} f \bigl( t_k^{(\alpha)} , p_k \bigr) 1_{ \{ f (
t_k^{(\alpha)} , p_k ) \le\ep\} } \lambda_{k+1}' .
\]
We can easily see that
$| Y_k | \le1$ and
$\sum_{k=n+1}^N Y_k$ is $\sF_N$-martingale.
By~\cite{Freedtail}, Theorem 1.6,
with~(\ref{eqreplace2}),
we obtain
\begin{eqnarray*}
&&\PP[ \alpha^2 ( \zt_\ep- n ) < \delta , \hat{\sg}_R  = \infty ]
\\
&&\qquad \le \PP\Biggl[ \sup_{n \le N \le N_\delta^{(\alpha)} } \sum _{k=n+1}^{N+1}
Y_k > \frac{\ep^2}{4 \alpha\sqrt{m+2}} \Biggr]
\\
&&\qquad \le \exp \biggl( - \frac{\ep^4} { 8 \sqrt{m+2} ( \alpha\ep^2 + 4 \alpha^2
\sqrt{m+2} ( N_\delta^{(\alpha)} - n ) ) } \biggr)
\\
&&\qquad \le \exp \biggl( - \frac{\ep^4} { 8 \sqrt{m+2} ( \alpha\ep^2 + 4 \sqrt
{m+2} \delta) } \biggr) .
\end{eqnarray*}
Hence~(\ref{eqtight}) follows
by combining this estimate with~(\ref{eqreplace}).
\end{pf*}
%

\section{Coupling by reflection}
\label{sec4}

For $k \in\R$,
let $U_{a,k}$ be a one-dimensional
Orn\-stein--Uhlenbeck process
defined as a solution to the following SDE:
\begin{eqnarray*}
d U_{a,k} (t) & = & - \frac{k}{2} U_{a,k} (t) \,dt + 2 \,d B (t) ,
\\
U_{a,k} ( T_1 ) & = & a.
\end{eqnarray*}
More explicitly,\vspace*{1pt}
$
U_{a,k} (t)
=
\mathrm{e}^{-k(t - T_1)/2} a + 2 \int_{T_1}^t \mathrm{e}^{k (s-t)/2} \,d B(s)
$.
Here $B(t)$ is the standard one-dimensional Brownian motion
as in the proof of Proposition~\ref{propnon-explosion}.
%
\begin{theorem} \label{thmain0}
Suppose
%
\begin{equation} \label{eqc-b0}
2 ( \nab Z (t) )^\flat+ \partial_t g (t) \le \Ric_{g(t)} + k g(t)
\end{equation}
holds for some $k \in\R$. Then, for each $x_1 , x_2 \in M$, there
exists a coupling $\mathbf{X} (t) := ( X_1 (t) , X_2 (t) )$ of two
$\sL_t$-diffusion processes starting at $(x_1 , x_2 )$ satisfying
\begin{eqnarray*} 
\PP\Bigl[ \inf_{T_1 \le t \le T} d_{g(t)} ( \mathbf{X}(t) ) > 0 \Bigr] & \le&
\PP\Bigl[ \inf_{T_1 \le t \le T} U_{d_{g(T_1)} (x_1, x_2),k} (t) > 0 \Bigr]
\\
& = & \chi \biggl( \frac{ d_{g(T_1)} (x_1 , x_2) } { 2 \sqrt{ \beta( T - T_1
) } } \biggr)
\end{eqnarray*}
for each $T \in[ T_1 , T_2 ]$,
where
\[
\chi(a) := \frac{1}{\sqrt{2 \pi}} \int_{-a}^a \mathrm{e}^{-u^2 /2}
\,du,\qquad
\beta(t) := \cases{ \displaystyle \frac{e^{kt} - 1}{k}, &\quad $k \neq0$,\vspace*{2pt}\cr
t, &\quad $k = 0$.}
\]
In addition, for $i = 1, 2$,
$X_i (t)$ is a solution to the martingale problem associated with
the time-inhomogeneous generator $\sL_t$ and
the filtration generated by $\mathbf{X}$.
\end{theorem}
\begin{remark}
\label{remO-U}
\textup{(i)}
Our assumption~(\ref{eqc-b0}) extends
existing curvature assumptions in two respects.
On the one hand,~(\ref{eqc-b0}) is nothing but~(\ref{eqc-b})
when $Z(t) \equiv0$ and $k = 0$.
On the other hand,~(\ref{eqc-b0}) can be regarded
as a natural extension of a lower Ricci curvature bound
by $k$.
Indeed, Bakry--\'{E}mery's curvature-dimension condition
$\mathsf{CD} (k ,\infty)$ (see~\cite{Bak97}, e.g.),
which is a natural extension of a lower Ricci curvature bound
by $k$, appears in~(\ref{eqc-b0})
when both~$Z(t)$ and $g(t)$ are independent of~$t$.\vadjust{\goodbreak}

\mbox{\hphantom{i}}\textup{(ii)}
Given $k > 0$,
a simple example satisfying~(\ref{eqc-b0}) can be constructed
by a scaling.
Indeed, for a complete metric $g$
whose Ricci curvature is nonnegative,
$g(t) = \mathrm{e}^{-k(t - T_1)} g$ satisfies~(\ref{eqc-b0})
when $Z(t) \equiv0$.

\textup{(iii)}
From the first item in this remark, when $Z(t) \equiv0$, one may
expect that~(\ref{eqc-b0}) works as an analog of Bakry--\'{E}mery's
$\mathsf{CD} (k , N)$ condition, which is equivalent to $\Ric_g \ge k$ and
$\dim M < N$ when $g(t)$ is independent of~$t$, instead of $\mathsf{CD}(k,
\infty)$ since $\dim M = m < \infty$ in our case. However, the
following observation suggests us that we should be careful: let us
consider~(\ref{eqc-b0}) in the case $k > 0$ and \mbox{$Z(t) \equiv0$}. When
$\partial_t g(t) \equiv0$, the Bonnet--Myers theorem tells us that the
diameter of $M$ is bounded and hence $M$ is compact. Moreover, the
Bonnet--Myers theorem still holds under $\mathsf{CD} (k, N)$
in the time-homogeneous case; see \cite
{Bak-LedBM-Sob,Bak-Qianvol,QianBonnetMyers}.
However, when $g(t)$ depends on~$t$,
it is no longer true
that~(\ref{eqc-b0}) implies the compactness of $M$.
In fact, we can easily obtain a noncompact $M$
enjoying~(\ref{eqc-b0}) with $k > 0$ for some $g(t)$
by following the observation
in the second item of this remark.
\end{remark}

By a standard argument,
Theorem~\ref{thmain0} implies the following estimate
for a gradient of the diffusion semigroup:
\begin{cor} \label{corsF}
Let $( ( X(t) )_{t \in[ T_1, T_2 ]} , ( \PP_x )_{x \in M} )$
be a $\sL_t$-diffusion process with $\PP_x [ X(T_1) = x ] = 1$.
For any bounded measurable function $f$ on $M$,
let us define $P_t f$ by $P_t f (x) := \E_x [ f ( X(t) ) ]$.
Then, under the same assumption as in Theorem~\ref{thmain0},
we have
\[
\limsup_{y \to x} \biggl| \frac{ P_t f (x) - P_t f (y) }{d_{g(T_1)} (x,y)} \biggr|
\le \frac{1}{\sqrt{2 \pi\beta( t - T_1 )}} \sup_{z , z' \in M} | f (z)
- f (z') | .
\]
In particular,
$P_t f$ is
$d_{g(T_1)}$-globally Lipschitz continuous
when $f$ is bounded.
\end{cor}
\begin{pf}
Let $\mathbf{X} = (X_1 , X_2 )$ be a coupling of $\sL_t$-diffusions $(
X(t) , \PP_x )$ and $( X(t) , \PP_y )$ given in Theorem~\ref{thmain0}.
Let $\tau^*$ be the coupling time of $\mathbf{X}$, that is, $ \tau^* :
= \inf\{ t \in[ T_1 , T_2 ] | \mathbf{X} (t) \in D (M) \} $. Let us
define $\mathbf{X}^* = ( X_1^* , X_2^* )$ of $( X(t) , \PP_x )$ and $(
X(t) , \PP_y )$ by
\[
\mathbf{X}^* (t) : =
\cases{\mathbf{X} (t), &\quad if $\tau^* > t$,\cr
( X_1 (t) , X_1 (t) ), &\quad otherwise.}
\]
Since $\tau^*$ is a stopping time
with respect to the filtration generated by
$\mathbf{X}$, and $X_i$ ($i = 1, 2$) is a solution
to the martingale problem associated with the
same filtration, $\mathbf{X}^*$ is again
a coupling of $\sL_t$-diffusion processes.
Since
$
\{ \tau^* > T \}
=
\{ \inf_{T_1 \le t \le T} d_{g(t)} ( \mathbf{X} (t) ) > 0 \}
$,
Theorem~\ref{thmain0} yields
\begin{eqnarray*}
P_t f (x) - P_t f (y)
& = & \E[ f ( X_1^* (t) ) - f ( X_2^* (t) ) ]
\\
& = & \E\bigl[ \bigl( f ( X_1^* (t) ) - f ( X_2^* (t) ) \bigr) 1_{ \{ \tau^* > t \} }
\bigr]
\\
& \le& \PP[ \tau^* > t ] \sup_{z,z' \in M} | f (z) - f(z') |
\\
& \le& \chi\biggl( \frac{d_{g(T_1)} (x,y)}{2 \sqrt{\beta( t - T_1 )}} \biggr)
\sup_{z,z' \in M} | f(z) - f(z') | .
\end{eqnarray*}
Hence the assertion holds
by dividing the both sides of the above inequality
by $d_{g(T_1)} (x,y)$
and
by letting $y \to x$ after that.
\end{pf}

%
As we did in the last section, let $( \gm_{xy}^{(t)} )_{x,y \in M}$ be
a measurable family of unit-speed minimal $g(t)$-geodesics such that
$\gm_{xy}^{(t)}$ joins $x$ and $y$. Without loss of generality, we may
assume that $\gm_{xy}^{(t)}$ is symmetric, that is, $ \gm_{xy}^{(t)} (
d_{g(t)} (x,y) - s ) = \gm_{yx}^{(t)} (s) $ holds.
Let us define
$\tilde{m}_{xy}^{(t)} \dvtx T_y M \to T_y M$
by
\[
\tilde{m}_{xy}^{(t)} v := v - 2 \bigl\langle v, \dot{\gm}_{xy}^{(t)}
\bigr\rangle_{g(t)} \dot{\gm}_{xy}^{(t)} \bigl( d_{g(t)} (x,y) \bigr) .
\]
This is a reflection with respect to a hyperplane
which is $g(t)$-perpendicular to $\dot{\gm}_{xy}^{(t)}$.
Let us define $m_{xy}^{(t)} \dvtx T_x M \to T_y M$ by
\[
m_{xy}^{(t)} (v) := \tilde{m}_{xy}^{(t)} \bigl( \bigl( \DD_{\gm_{xy}^{(t)}}^{(t)}
v \bigr) \bigl( d_{g(t)} (x,y) \bigr) \bigr) .
\]
Clearly $m_{xy}^{(t)}$ is a $g(t)$-isometry.
As in the last section,
let $\Ph^{(t)} \dvtx M \to\mathscr{O}^{(t)} (M)$ be
a measurable section of
the $g(t)$-orthonormal frame bundle
$\mathscr{O}^{(t)} (M)$ of~$M$.
Let us define two measurable maps
$\Phi_i^{(t)} \dvtx M \times M \to\mathscr{O}^{(t)} (M)$
for $i=1,2$ by
\begin{eqnarray*}
\Ph_1^{(t)} (x,y) &
:= &\Ph^{(t)} (x),
\\
\Ph_2^{(t)} (x,y)
& := &
\cases{
m_{xy}^{(t)} \Ph_1^{(t)} (x,y),&\quad
$(x,y) \in M \times M \setminus D(M)$,
\vspace*{2pt}\cr
\Phi^{(t)} (x), &\quad $(x,y) \in D (M)$.}
\end{eqnarray*}
Take $x_1 , x_2 \in M$.
By using $\Phi_i^{(t)}$,
we define a coupled geodesic random walk
$\mathbf{X}^\alpha(t) = ( X_1^\alpha(t) , X_2^\alpha(t) )$
by $X^\alpha_i (T_1) = x_i$ and,
for $t \in[ t_n^{(\alpha)} , t_{n+1}^{(\alpha)} ]$,
\begin{eqnarray*}
\tilde{\xi}_{n+1}^i & : = & \sqrt{m+2} \Phi_i^{( t_n^{(\alpha)} )} \bigl(
\mathbf{X}^\alpha\bigl( t_n^{(\alpha)} \bigr) \bigr) \xi_{ n + 1 } ,
\\ 
X_i^\alpha( t ) & := & \exp_{X_i^\alpha( t_n^{(\alpha)} )}^{(
t_n^{(\alpha)} )} \biggl( \frac{ t - t_n^{(\alpha)} }{\alpha^2} \bigl(
\alpha\tilde{\xi}_{n+1}^i + \alpha^2 Z\bigl( t_n^{(\alpha)} \bigr) \bigr) \biggr)
\end{eqnarray*}
for $i = 1, 2$.
We can easily verify that $X_i^\alpha$ has
the same law as $X^\alpha$ with $x_0 = x_i$.
%
%

In what follows, we assume~(\ref{eqc-b0}).
We can easily verify that
it implies Assumption~\ref{assnon-explosion}.
Thus, by Theorem~\ref{thIP},
$( \mathbf{X}^\alpha)_{\alpha> 0}$ is tight
under Assumption~\ref{assnon-explosion}.
In addition, a subsequential limit
$\mathbf{X}^{\alpha_k} \to\mathbf{X} = ( X_1 , X_2 )$
in law exists, and it is a coupling of
two $\sL_t$-diffusion processes
starting at $x_1$ and $x_2$, respectively.
We fix such a subsequence\vadjust{\goodbreak} $( \alpha_k )_{k \in\N}$.
In the rest of this paper,
we use the same symbol $\mathbf{X}^\alpha$
for the subsequence $\mathbf{X}^{\alpha_k}$
and the term ``$\alpha\to0$'' always means
the subsequential limit ``$\alpha_k \to0$.''

We will prove that the coupling $\mathbf{X}$ obtained as above is a
desired one in Theorem~\ref{thmain0}. We first remark that we can
easily verify that $X_i$ ($i = 1, 2$) is a solution to the martingale
problem associated with the filtration generated by $\mathbf{X}$ in the
same way as in the proof of Theorem~\ref{thIP}. Set $ \hat{\sg}_R^{i} :
= \sg_R ( d_{g(\cdot)} ( o , X_i^\alpha(\cdot) ) ) $ for $i = 1,2$. We
fix $R > 1$ sufficiently large until the beginning of the proof of
Theorem~\ref{thmain0}. Let $M_0 \subset M$ be a relatively compact open
set satisfying~(\ref{eqbddcpt}) for $2R$ instead of $R$. We next show a
difference inequality of $d_{g(t)} ( \mathbf{X}^\alpha(t) )$. To
describe\vspace*{-1pt} it, we will introduce several notation as in the last section.
For simplicity, let us denote $\gm_{X_1^\alpha( t_n^{(\alpha)} )
X_2^\alpha( t_n^{(\alpha)} )}^{( t_n^{(\alpha)} )}$ by $\bar{\gm}_n$.
Let us define a vector field $V_{n+1}$ along $\bar{\gm}_n$ by
\[
V_{n+1} : = \DD_{\bar{\gm}_n}^{(t_n^{(\alpha)})} \bigl( \tilde{\xi}_{n+1}^1
- \langle\tilde{\xi}_{n+1}^1 , \dot {\bar{\gm}}_n \rangle_{g
(t_n^{(\alpha)})} \dot{\bar{\gm}}_n (0) \bigr) .
\]
Take $v \in\R^m$.
Let us define $\lambda_{n+1}^*$ and $\Lm_{n+1}^*$ by
\begin{eqnarray*}
\lambda_{n+1}^*
& : = & \cases{ 2 \langle\tilde{\xi}_{n+1}^1 ,
\dot{\bar{\gm}}_n \rangle _{g(t_n^{(\alpha)})}, &\quad if $(y_1 , y_2)
\notin D(M)$,\vspace*{2pt}\cr 2 \sqrt{m+2} \langle\xi_{n+1} , v \rangle, &\quad
otherwise,}
\\
\Lm_{n+1}^*
& : = &
\frac12 \biggl(
\int_0^{d_{g(t_n^{(\alpha)})} ( \mathbf{X}^\alpha(t_n^{(\alpha)}) )}
\bigl( \partial_t g \bigl( t_n^{(\alpha)} \bigr) + 2 \bigl( \nab Z \bigl( t_n^{(\alpha
)} \bigr) \bigr)^\flat\bigr)
\\
&&\hspace*{83.8pt}{}\times( \dot{\bar{\gm}}_n (s) , \dot{\bar{\gm}}_n (s) )
\,ds
\\
&&\hspace*{135pt}{} +
I_{\bar{\gm}_n}^{(t_n^{(\alpha)})}
( V_{n+1} )
\biggr) 1_{\{ \mathbf{X}^\alpha( t_n^{(\alpha)} ) \notin D(M) \} }.
\end{eqnarray*}
For $\delta\ge0$, let us define
$\tau_\delta \dvtx \sC_1 \to[ T_1 , T_2 ] \cup\{ \infty\}$
by
\[
\tau_\delta(w) := \inf\{ t \ge T_1 | w(t) \le\delta \}.
\]
We also define $ \hat{\tau}_\delta $ by $ \hat{\tau}_\delta : =
\tau_\delta( d_{g(\cdot)} ( \mathbf{X}^\alpha(\cdot) ) ) . $
%
\begin{lem} \label{lemc-2var}
For $n \in\N_0$ with $n < N^{(\alpha)}$,
we have
%
\begin{eqnarray}\label{eqc-2var}
\mathrm{e}^{ k t_{n+1}^{(\alpha)} / 2}
d_{g( t_{n+1}^{(\alpha)} )} \bigl( \mathbf{X}^\alpha\bigl( t_{n+1}^{(\alpha)}
\bigr) \bigr)
& \le&
\biggl( 1 + \frac{k}{2} \biggr)
\mathrm{e}^{ k t_n^{(\alpha)} / 2}
d_{g( t_n^{(\alpha)} )} \bigl( \mathbf{X}^\alpha\bigl( t_n^{(\alpha)} \bigr) \bigr)
\nonumber\\
&&{} +
\mathrm{e}^{ k t_n^{(\alpha)} / 2} ( \alpha\lambda_{n+1}^* + \alpha
^2 \Lm
_{n+1}^* )\\
&&{} + o (\alpha^2 ),\nonumber
\end{eqnarray}
when
$n < \hat{\tau}_\delta\wg\hat{\sg}_R^1 \wg\hat{\sg}_R^2$
and $\alpha$ is sufficiently small.
Moreover, we can control the error term $o(\alpha^2)$
uniformly in the position of $\mathbf{X}^\alpha$.
\end{lem}
\begin{pf}
When $( t_n^{(\alpha)} , \mathbf{X}^\alpha( t_n^{(\alpha)} ) )
\notin\Cutst$,
(\ref{eqc-2var}) is just a consequence of
the second variational formula for the distance function\vadjust{\goodbreak}
combined with the index lemma for $I_{\bar{\gm}_n}^{(t_n^{(\alpha)})}$.
To include the case
$( t_n^{(\alpha)} , \mathbf{X}^\alpha( t_n^{(\alpha)} ) ) \in
\Cutst$ and
to obtain a uniform control of $o (\alpha^2)$,
we extend this argument.
Let us define~$H$ and
$p_1 , p_2 \dvtx H \to[ T_1 , T_2 ] \times\bar{M}_0 \times\bar{M}_0$
by
\begin{eqnarray*}
&
\displaystyle H : = \bigl\{ (t,x,y,z) |
t \in[ T_1 , T_2 ], x,y,z \in\bar{M}_0 ,d_{g(t)} (x,y) \ge\delta,&\\[-2pt]
&\hspace*{46.6pt}\hspace*{46.6pt} d_{g(t)} (x,y) = 2 d_{g(t)} (x,z) = 2
d_{g(t)} (y,z)\bigr\},&
\\[-2pt]
&p_1 (t,x,y,z)
: =
(t,x,z),&
\\[-2pt]
&p_2  (t,x,y,z)
: =
(t,y,z).&
\end{eqnarray*}
If $\mathbf{q} = (t,x,y,z) \in H$,
then
$p_1 (\mathbf{q} ) , p_2 ( \mathbf{q} ) \notin\Cutst$
since $z$ is on a midpoint of
a minimal $g(t)$-geodesic joining $x$ and $y$.
Since $H$ is compact,
$p_1 (H)$ and $p_2 (H)$ are also compact.
Hence there is a constant $\h> 0$ such that
\begin{eqnarray*}
&&\inf \bigl\{ | t - t' | + d_{g(t)} ( x, x' ) + d_{g(t)} ( y, y' ) | (t,x,y)
\in p_1 (H) \cup p_2 (H) ,\\[-2pt]
&&\hspace*{212pt} (t',x',y') \in \Cutst\bigr\}> \h.
\end{eqnarray*}
Take $\alpha> 0$ sufficiently small
relative to $\h$ and $\delta$.
Set
\begin{eqnarray*}
p_n
& := &
\bar{\gm}_n \biggl( \frac{d_{g ( t_n^{(\alpha)} )} ( \mathbf
{X}^\alpha( t_n^{(\alpha)} ) ) }{2} \biggr),
\\[-2pt]
p_n '
& := &
\exp_{p_n}^{( t_n^{(\alpha)} )}
\biggl( V_{n+1} \biggl( \frac{d_{g( t_n^{(\alpha)} ) } ( \mathbf
{X}^\alpha( t_n^{(\alpha)} ) ) }{2} \biggr) \biggr)
.
\end{eqnarray*}
By the triangle inequality,
we have
\begin{eqnarray*}
d_{g ( t_{n}^{(\alpha)} )} \bigl( \mathbf{X}^\alpha\bigl( t_{n}^{(\alpha)} \bigr) \bigr)
& = &
d_{g( t_{n}^{(\alpha)} )}
\bigl( X_1^\alpha\bigl( t_{n}^{(\alpha)} \bigr) , p_n \bigr)
+
d_{g( t_{n}^{(\alpha)} )}
\bigl( p_n , X_2^\alpha\bigl( t_{n}^{(\alpha)} \bigr) \bigr) ,
\\[-2pt]
d_{g ( t_{n+1}^{(\alpha)} )} \bigl( \mathbf{X}^\alpha\bigl( t_{n+1}^{(\alpha
)} \bigr) \bigr)
& \le&
d_{g( t_{n+1}^{(\alpha)} )}
\bigl( X_1^\alpha\bigl( t_{n+1}^{(\alpha)} \bigr) , p_n' \bigr)
+
d_{g( t_{n+1}^{(\alpha)} )}
\bigl( p_n' , X_2^\alpha\bigl( t_{n+1}^{(\alpha)} \bigr) \bigr).
\end{eqnarray*}
Since
$
(
t_n^{(\alpha)} ,
\mathbf{X}^\alpha( t_n^{(\alpha)} ),
p_n
)
\in H
$,
we can apply the second variation formula
to each term on the right-hand side of
the above inequality.
Hence we obtain~(\ref{eqc-2var}).
For a uniform control of the error term,
we remark that $\bar{\gm}_n$ is included in $M_0$,
and the $g(t_n^{(\alpha)})$-length of $\bar{\gm}_n$ is
bigger than $\delta$.
These facts follows from
$n < \hat{\tau}_\delta\wg\hat{\sg}_R^1 \wg\hat{\sg}_R^2$
and the choice of $M_0$.
Thus the every calculation of the second variation
formula above is done
on a compact subset of $[ T_1 , T_2 ] \times M_0 \times M_0$
which is uniformly away from $\Cutst$.
It yields the desired result.\vspace*{-2pt}
\end{pf}

Let us define
a continuous stochastic process $U_a^\alpha$
on $\R$ starting at $a$
by
\[
U_a^\alpha(t) := \mathrm{e}^{-kt/2} a + \alpha\mathrm{e}^{-kt/2} \Biggl(
\sum_{j=1}^{n} \mathrm{e}^{k t_j^{(\alpha)} /2} \lambda_j^* + \frac{ t
- t_n^{(\alpha)} }{\alpha^2} \mathrm{e}^{ k t_n^{(\alpha )} /2}
\lambda_{n + 1}^* \Biggr) .
\]
As a final preparation of the proof of Theorem~\ref{thmain0},
we show the following comparison theorem
for the distance process of coupled geodesic random walks.\vadjust{\goodbreak}
\begin{lem} \label{lemchain-LLN}
For each $\ep> 0$, there exists a family of events
$( E_\ep^\alpha)_{\alpha}$ such that
$
\PP[ E_\ep^\alpha]
$
converges to 1
as $\alpha\to0$
and
%
\begin{equation} \label{eqchain}
d_{g(t)} ( \mathbf{X}^\alpha(t) ) \le U^\alpha_{d_{g(T_1)} (
\mathbf{X}^\alpha(T_1) )} (t) + \ep
\end{equation}
for all $ t \in [ T_1 , T_2 \wg \hat{\tau}_\delta\wg \hat{\sg}_R^1 \wg
\hat{\sg}_R^2 ] $ on $E_\ep^\alpha$ for sufficiently small $\alpha$.
\end{lem}
\begin{pf}
In a similar way as in the proof of Lemma~\ref{lemcomparison1}, we can
complete the proof once we have found $E_\ep^\alpha$ on which
(\ref{eqchain}) holds when $ t = t_n^{(\alpha)} \in [ T_1 , T_2 \wg
\hat{\tau}_\delta\wg \hat{\sg}_R^1 \wg \hat{\sg}_R^2 ] $. Set $
\bar{\Lm}_{n+1}^* := \E[ \Lm_{n+1}^* | \sF_n ] $.
Then $\sum_{j=1}^n \mathrm{e}^{ k t_{j-1}^{(\alpha)} /2 } ( \Lm_j^* -
\bar{\Lm}_j^* )$ is
an $\sF_n$-local martingale.
Indeed, $\Lm_{n+1}^*$ is bounded
if $t_n^{(\alpha)} < \hat{\sg}_R^1 \wg\hat{\sg}_R^2$ and
so is $\bar{\Lm}_{n+1}^*$.
Let us define $E_\ep^\alpha$ by
\[
E_\ep^\alpha : = \Biggl\{ \mathop{\sup_{N \le N^{(\alpha)}}}_{t_N^{(\alpha)}
\le T_2 \wg\hat{\sg}_R^1 \wg\hat{\sg}_R^2} \sum_{j=1}^{N+1} \mathrm
{e}^{k t_j^{(\alpha)} / 2} ( \Lm_j^* - \bar{\Lm}_j^* ) \le\frac{\ep}{2
\alpha^2 ( T_2 - T_1 )} \Biggr\} .
\]
In a similar way as in Lemma~\ref{lemLLN} or~\cite{K8}, Lemma 6,
$\lim_{\alpha\to0} \PP[ E_\ep^\alpha] = 1$ holds. On~$E_\ep^\alpha$, we
can replace $\alpha^2 \mathrm{e}^{k t^{(\alpha)}_n / 2} \Lm_{n+1}^*$ in
(\ref{eqc-2var}) with $\alpha^2 \mathrm{e}^{k t^{(\alpha)}_n / 2}
\bar{\Lm}_{n+1}^* + \ep/ ( 2(T_2 - T_1) )$. Since we have $ (m+2) \E [
\langle\xi_i , e_k \rangle \langle\xi_i , e_l \rangle ] = \delta_{kl}
$, we obtain
\[
\bar{\Lm}_{n+1}^* \le- \frac{k}{2} d_{g( t_n^{(\alpha)} )} \bigl(
\mathbf{X}^\alpha\bigl( t_n^{(\alpha)} \bigr)\bigr).
\]
Thus an iteration of Lemma~\ref{lemc-2var} implies
(\ref{eqchain}) on $E_\ep^\alpha$
when $t = t_n^{(\alpha)}$.
\end{pf}
%
%
\begin{pf*}{Proof of Theorem~\ref{thmain0}}
Take $\ep\in( 0 , 1 )$ arbitrarily.
Let $R > 1$ be sufficiently large
so that
\[
\limsup_{\alpha\to0} \PP[ \hat{\sg}_R^1 \wg\hat{\sg}_R^2 \le T_2 ] <
\ep .
\]
It is possible by Proposition~\ref{propnon-explosion}.
Set $a := d_{g(T_1)} ( x_1, x_2 )$.
Take $T \in[ T_1 , T_2 ]$, and
let $\delta> 0$ be $\delta> 2 \ep$.
Then Lemma~\ref{lemchain-LLN} yields
\begin{eqnarray*}
\PP[ \hat{\tau}_\delta> T ] & \le& \PP[ \{ \hat{\tau}_\delta> T \} \cap
E_\ep ^\alpha\cap\{ \hat{\sg}^1_R \wg\hat{\sg}_R^2 > T \} ] + 2 \ep
\\
& \le& \PP[ \tau_{\delta/ 2} ( U^\alpha_a ) > T ] + 2 \ep .
\end{eqnarray*}
Thus we obtain
\[
\limsup_{\alpha\to0} \PP[ \hat{\tau}_\delta> T ] \le
\limsup_{\alpha\to0} \PP\Bigl[ \inf_{t\in[ T_1 ,T ]} U^\alpha_a (t)
\ge\delta/ 2 \Bigr]
\]
by letting $\ep\downarrow0$.
Note that $U^\alpha_a$ converges in law
to $U_a$ as $\alpha\to0$.
Since
\[
\bigl\{ \mathbf{w} \in C ( [ T_1, T_2 ] \to M \times M ) |
\tau_\delta\bigl(d_{g(\cdot)} ( \mathbf{w} (\cdot) ) \bigr) > T \bigr\}
\]
is open,
and
$\{ w | \inf_{t\in[ T_1 ,T_2 ]} w(t) \ge\delta/ 2
\}$
is closed in $C ( [ 0, T ] \to\R)$,
the Portmanteau theorem yields
\begin{eqnarray*}
\PP\Bigl[ \inf_{T_1 \le t \le T} d_{g(t)}( \mathbf{X}( t ) ) > \delta\Bigr] &
\le& \liminf_{\alpha\to0} \PP[ \hat{\tau}_\delta> T ]
\\
& \le& \limsup_{\alpha\to0} \PP\Bigl[ \inf_{t\in[ T_1 ,T ]} U^\alpha_a (t)
\ge\delta/ 2 \Bigr]
\\
& \le& \PP\Bigl[ \inf_{t\in[ T_1 ,T ]} U_a (t) \ge\delta/ 2 \Bigr] .
\end{eqnarray*}
Therefore the conclusion follows
by letting $\delta\downarrow0$.
\end{pf*}

\section{Coupling by parallel transport}\label{sec5}

As a final part of the paper,
we will see that we can also construct a coupling
by parallel transport by following our manner.
In the construction of the coupling by reflection,
we used a~map~$m_{xy}^{(t)}$.
By following the same argument
after omitting $\tilde{m}_{xy}^{(t)}$
in the definition of $m_{xy}^{(t)}$,
we obtain a coupling by parallel transport.
The difference of it
from the coupling by reflection
is the absence of the term corresponding to $\lambda_n^*$,
which comes from the first variation of arc length.
As a result, we can show the following;
cf.~\cite{K8}:
\begin{theorem} \label{thparallel}
Assume~(\ref{eqc-b0}).
For $x_1 , x_2 \in M$, there is
a coupling $\mathbf{X}(t) = ( X_1 (t) , X_2 (t) )$
of two $\sL_t$-diffusion processes
starting at $x_1$ and $x_2$ at time~$T_1$, respectively,
such that
\[
d_{g(t)} ( \mathbf{X} (t) ) \le \mathrm{e}^{-k(t-s)/2} d_{g(s)} (
\mathbf{X} (s) )
\]
for $T_1 \le s \le t \le T_2$ almost surely.
\end{theorem}

It recovers a part of results
studied in~\cite{Arn-Coul-Thalhoriz}.
In particular,
a contraction type estimate
for Wasserstein distances
under the heat flow follows.
\begin{pf*}{Proof of Theorem~\ref{thparallel}}
Let us construct a coupling
by parallel transport of geodesic random walks
$
\mathbf{X}^\alpha
=
( X_1^\alpha, X_2^\alpha)
$
starting at $(x_1 , x_2 ) \in M \times M$
by following the procedure
stated just before Theorem~\ref{thparallel}.
By taking a subsequence,
we may assume that
$\mathbf{X}^\alpha$ converges in law as $\alpha\to0$.
We denote the limit by $\mathbf{X} = ( X_1 , X_2 )$.
In what follows, we prove
\[
\PP\Bigl[ \sup_{T_1 \le s \le t \le T_2} \bigl( \mathrm{e}^{kt/2} d_{g(t)} (
\mathbf{X} (t) ) - \mathrm{e}^{ks/2} d_{g(s)} ( \mathbf {X} (s) ) \bigr) >
\ep\Bigr] = 0
\]
for any $\ep> 0$.
By virtue of the Portmanteau theorem
together with Proposition~\ref{propnon-explosion},
it suffices to show
%
\begin{eqnarray} \label{eqpara0}
&&\lim_{\alpha\to0} \PP\Bigl[ \sup_{T_1 \le s \le t \le T_2} \bigl(
\mathrm{e}^{kt/2} d_{g(t)} ( \mathbf{X}^\alpha(t) ) - \mathrm{e}^{ks/2}
d_{g(s)} ( \mathbf{X}^\alpha(s) ) \bigr)
> \ep,
\nonumber\\[-8pt]\\[-8pt]
&&\hspace*{208.2pt}
\hat{\sg}_R^1 \wg\hat{\sg}_R^2 = \infty \Bigr] = 0 \nonumber
\end{eqnarray}
for any $R > 1$. We write $ d_n := \mathrm{e}^{k t_n^{(\alpha)} / 2}
d_{g(t_n^{(\alpha)})} ( \mathbf {X}^\alpha (t_n^{(\alpha)}) ) $ in this
proof for simplicity of notation. For $\delta> 0$, let us define a
sequence of $\sF_n$-stopping times~$S_l$ by $S_0 : = 0$ and
\begin{eqnarray*}
S_{2l+1} & : = & \inf \{ j \ge S_{2l} | d_j \le\delta\} \wg
N^{(\alpha)} ,
\\
S_{2l} & : = & \inf \{ j \ge S_{2l-1} | d_j \ge2 \delta \} \wg
N^{(\alpha)} .
\end{eqnarray*}
Note that $d_{S_{2l-1}} \le3\delta$ holds
on $\{ \hat{\sg}_R^1 \wg\hat{\sg}_R^2 = \infty\}$
for sufficiently small $\alpha$.
As mentioned
just before Theorem~\ref{thparallel},
the difference inequality~(\ref{eqc-2var}) holds with $\lambda^* = 0$
when
$
S_{2l-1}
\le
n
<
S_{2l}
\wg\bar{\sg}_R^1
\wg\bar{\sg}_R^2
$
for some $l \in\N_0$.
In this case, the error term $o( \alpha^2 )$ is
controlled uniformly also in $l$.
Let us define an event $E^\alpha_\delta$ by
\[
E_\delta^\alpha : = \Biggl\{ \mathop{\sup_{n \le N \le
N^{(\alpha)}}}_{t_N^{(\alpha)} \le T_2 \wg\hat{\sg}_R^1
\wg\hat{\sg}_R^2} \sum_{j=n+1}^{N+1} \mathrm{e}^{k t_j^{(\alpha)} / 2}
( \Lm_j^* - \bar{\Lm}_j^* ) \le\frac{\delta}{2 \alpha^2} \Biggr\}.
\]
Then, as in Lemmas~\ref{lemLLN} and~\ref{lemchain-LLN},
we can show $\lim_{\alpha\to0} \PP[ E_\delta^\alpha] = 1$.
On $E_\delta^\alpha\cap\{ \hat{\sg}_R^1 \wg\hat{\sg}_R^2 =
\infty\}$,
we have $d_N \le d_n + \delta$
for $S_{2l-1} \le n \le N \le S_{2l}$
if $\alpha$ is sufficiently small.
Moreover, for $n < S_{2l-1} \le N < S_{2l}$,
\[
d_N - d_n \le ( d_N - d_{S_{2l-1}} ) + d_{S_{2l-1}} \le 5 \delta.
\]
In the case $S_{2l} \le N < S_{2l+1}$,
we obtain $d_N - d_n \le2 \delta$.
Thus $d_N - d_n \le5 \delta$ holds for all $n < N$
on $E_\delta^\alpha\cap\{ \hat{\sg}_R^1 \wg\hat{\sg}_R^2 =
\infty\}$.
Take $\delta> 0$ less than $\ep/ 10$.
Then our observations yield~(\ref{eqpara0})
since
$
d_{g(t)} ( \mathbf{X}^\alpha(t) )
-
d_{g( \lfloor t \rfloor_\alpha)} ( \mathbf{X} ( \lfloor t \rfloor
_\alpha) )
$
becomes uniformly small
on $\{ \hat{\sg}_R^1 \wg\hat{\sg}_R^2 = \infty\}$
as $\alpha\to0$.
\end{pf*}



%
\printaddresses

\end{document}